\documentclass[12pt,twoside,reqno]{amspaper}
\pdfoutput=1 
\usepackage[edges]{mathdef}
\usepackage{macros}
\usepackage{natbib}
\bibpunct{(}{)}{,}{a}{}{,}

\newcommand{\yclass}{\mathscr{Y}}
\renewcommand{\pastx}{\mathbf{X}}
\renewcommand{\ffilt}[1]{\mathscr{X}_{\!{#1}}}
\usepackage{comment}
\begin{document}
\title{Graphical modelling of multivariate time series}
\author{Michael Eichler}
\thanks{{\em E-mail address:} m.eichler@maastrichtuniversity.nl (M.~Eichler)}
\affil{Department of Quantitative Economics,
Maastricht University\\P.O.~Box 616, 6200 MD Maastricht, The Netherlands}
\dedicatory{{\upshape\today}}
\begin{abstract}
We introduce graphical time series models for the analysis
of dynamic relationships among variables in multivariate
time series. The modelling approach is based on the notion of strong
Granger causality and can be applied to time series with non-linear
dependencies. The models are derived from ordinary time series models
by imposing constraints that are encoded by mixed graphs.
In these graphs each component series is represented
by a single vertex and directed edges indicate possible Granger-causal
relationships between variables while undirected edges are used to
map the contemporaneous dependence structure. We introduce
various notions of Granger-causal Markov properties and
discuss the relationships among them and to other Markov
properties that can be applied in this context. Examples
for graphical time series models include nonlinear autoregressive
models and multivariate ARCH models.
\\[1em]
\end{abstract}
\maketitle

\section{Introduction}

Graphical models have become an important tool for the statistical
analysis of complex multivariate data sets, which are now increasingly
available in many scientific fields. The key feature of these models
is to merge the probabilistic concept of conditional independence with
graph theory by representing possible dependences among the variables
of a multivariate distribution in a graph. This has led to simple
graphical criteria for identifying the conditional independence relations
that are implied by a model associated with a given graph.
Further important advantages of the graphical modelling approach are
statistical efficiency due to parsimonious parameterizations of the joint
distribution of the variables and the visualization of complex dependence
structures, which allows an intuitive understanding of the interrelations
among the variables and, thus, facilitates the communication of statistical
results. For an introduction to graphical models we refer to the
monographs by \citet{whittaker90}, \citet{edwardsgraphs}, and
\citet{coxwermuth96}; a mathematically more rigorous treatment can be found
in \citet{SL96}.

While graphical models originally have been developed for variables that
are sampled with independent replications, they have been applied more
recently also to the analysis of time dependent data. Some first general
remarks concerning the potential use of graphical models in time series
analysis can be found in \citet{brill96}; since then there has been an
increasing interest in the use of graphical modelling techniques
for analyzing multivariate time series \citep[e.g.,][]{stanghellini99,RD00,
reale01,hsss,oxleyreale04,moneta05,eichlerhandbook,eichlerpathdiagr}.
However, all these works have been restricted to the analysis of linear
interdependences among the variables whereas the recent trend in time
series analysis has shifted towards non-linear parametric and non-parametric
models \citep[e.g.,][]{tong93,rothman99,fanyao03}. Moreover, in most of these
approaches, the variables at different time points are represented by
separate nodes, which leads to graphs with theoretically infinitely
many vertices for which no rigorous theory exists so far.

In this paper, we present a general approach for graphical modelling
of multivariate stationary time series, which is based on simple
graphical representations of the dynamic dependences of a process.
To this end, we utilize the concept of strong Granger causality
\citep[e.g.,][]{fm82}, which is formulated in terms of conditional
independences and, thus, can be applied to model arbitrary non-linear
relationships among the variables. The concept of Granger causality
originally has been introduced by \citet{granger69} and is commonly
used for studying dynamic relationships among the variables in
multivariate time series.

For the graphical representations, we
consider mixed graphs in which each variable as a complete time series
is represented by a single vertex and directed edges indicate possible
Granger-causal relationships among the variables while undirected edges
are used to map the contemporaneous dependence structure.
We note that similar graphs have been used in \citet{eichlerpathdiagr}
as path diagrams for the autoregressive structure of weakly stationary
processes or---without undirected edges---in \citet{didelez07sjs} for
graphical modelling of time-continuous composable finite Markov processes
based on the concept of local independence \citep{aalen87}. Formally, the
graphical encoding of the dynamic structure of
a time series is achieved by a new type of Markov properties, which we
call Granger-causal Markov properties. We introduce various levels,
namely the pairwise, the local, the block-recursive, and the global
Granger-causal Markov property, and discuss the relationships among them.
In particular, we give sufficient conditions under which
the various Granger-causal Markov properties are equivalent;
such conditions allow formulating models based on a simple Markov property
while interpreting the associated graph by use of the global
Granger-causal Markov property.

The paper is organized as follows. In Section \ref{sect:graphmodel},
we introduce the concepts of Granger-causal Markov properties and
graphical time series models; some examples of graphical time series
models are presented in Section \ref{sect:examples}. In Section
\ref{sect:globalMP}, we discuss global Markov properties, which relate
certain separation properties of the graph to conditional independence
or Granger noncausality relations among the variables of the process.
Finally in Section \ref{sect:discussion}, we compare the presented
graphical modelling approach with other approaches in the literature and
discuss possible extensions. The proofs are technical and put into the
appendix.

\section{Graphical time series models}
\label{sect:graphmodel}
In graphical modelling, the focus is on multivariate statistical models for
which the possible dependences between the studied variables can be
represented by a graph. In multivariate time series analysis, statistical models for a time series $X_V=\big(X_V(t)\big)_{t\in\znum}$ are usually specified in terms of the conditional distribution of $X_V(t+1)$ given its past
$\pastx_V(t)=\big(X_V(s)\big)_{s\leq t}$ in order to study the dynamic
relationships over time among the series. Thus, a time series model may be
described formally as a family of probability kernels $P$ from
$\rnum^{V\times\nnum}$ to $\rnum^V$, and we write $X_V\sim P$ if $P$ is a
version of the conditional probability of $X_V(t+1)$ given $\pastx_V(t)$.

For modelling specific dependence structures, we utilize the concept
of Granger (non-)causality, which has been introduced by \citet{granger69}
and has proved to be particularly useful for studying dynamic relationships
in multivariate time series. This probabilistic concept of noncausality
from a process $X_a$ to another process $X_b$ is based on studying
whether at time $t$ the next value of $X_b$ can be better predicted
by using the entire information up to time $t$ than by using the same
information apart from the former series $X_a$. In practice, not all
relevant variables may be available and, thus, the notion of
Granger causality clearly depends on the used information set.
In the sequel, we use the concept of strong Granger noncausality
\citep[e.g.,][]{fm82}, which is defined in terms of conditional independence
and $\sigma$-algebras and, thus, can be used also for non-linear time
series models.

Let $X_V=\big(X_V(t)\big)_{t\in\znum}$ with
$X_V(t)=(X_v(t))_{v\in V}\in\rnum^{V}$ be a multivariate stationary
stochastic process on a probability space $(\Omega,\falg,\prob)$.
For $A\subseteq V$, we denote by $X_A=(X_A(t))_{t\in\znum}$ the
multivariate subprocess with components $X_a$, $a\in A$. The information
provided by the past and present values of $X_A$ at time $t\in\znum$ can
be represented by the sub-$\sigma$-algebra $\ffilt{A}(t)$ of $\falg$
that is generated by $\pastx_A(t)=\big(X_A(s)\big)_{s\leq t}$. We write
$\ffilt{A}=(\ffilt{A}(t),t\in\znum)$ for the filtration induced by $X_A$.
This leads to the following definition of strong Granger noncausality
in multivariate time series; for ease of notation, we subsequently usually
drop the attribute ``strong''.

\begin{definition}
\label{def-causality}
Let $A$ and $B$ be disjoint subsets of $V$.
\begin{romanlist}
\item
$X_A$ is {\em strongly Granger-noncausal} for $X_B$ with respect
to the filtration $\ffilt{V}$ if
\[
\ffilt{B}(t+1)\indep\ffilt{A}(t)\given\ffilt{V\without A}(t)
\]
for all $t\in\znum$. This will be denoted by
$X_A\noncausal X_B\wrt{\ffilt{V}}$.
\item
$X_A$ and $X_B$ are {\em contemporaneously conditionally independent}
with respect to the filtration $\ffilt{V}$ if
\[
\ffilt{A}(t+1)\indep\ffilt{B}(t+1)\given
\ffilt{V}(t)\vee\ffilt{V\without(A\cup B)}(t+1)
\]
for all $t\in\znum$. This will be denoted by $X_A\noncorr X_B\wrt{\ffilt{V}}$.
\end{romanlist}
\end{definition}

Intuitively, the dynamic relationships of a stationary multivariate time series $X_V$
can be visualized by a mixed graph $G=(V,E)$ in which each vertex $v\in V$
represents one component $X_v$ and two vertices $a$ and $b$ are joined by
a directed edge $a\DE b$ whenever $X_a$ is Granger-causal for $X_b$ or
by an undirected edge $a\UE b$ whenever $X_a$ and $X_b$ are contemporaneously
conditionally dependent. Conversely, for formulating models
with specific dynamic dependences, a mixed graph $G$ can be associated
with a set of Granger noncausality and contemporaneous
conditional independence constraints that are imposed on a time series
model for $X_V$. Such a set of conditional independence relations encoded
by a graph $G$ is generally known as Markov property with respect to $G$.
In the context of multivariate time series, graphs may encode different
types of conditional independence relations, and we therefore speak of
Granger-causal Markov properties when dealing with Granger noncausality
and contemporaneous conditional independence relations. In the following
definition, $\parent{a}=\{v\in V|v\DE a\in E\}$ denotes the set of
parents of a vertex $a$, while $\neighbour{a}=\{v\in V|v\UE a\in E\}$
is the set of neighbours of $a$; furthermore, for $A\subseteq V$, we define
$\parent{A}=\cup_{a\in A}\parent{a}\without A$ and
$\neighbour{A}=\cup_{a\in A}\neighbour{a}\without A$ .

\begin{definition}[Granger-causal Markov properties]
\label{markovprop}
Let $G=(V,E)$ be a mixed graph. Then the stochastic process $X_V$ satisfies
\begin{list}{{}}{\usecounter{alphcount}
\labelwidth2.3em
\leftmargin2.8em\labelsep0.5em\topsep0.25em plus 0.5ex
\itemsep0.25em plus 0.5ex\parsep0em}
\item[(PC)]
the {\em pairwise Granger-causal Markov property}
with respect to $G$ if for all $a,b\in V$ with $a\neq b$
\begin{romanlist}
\item
$a\DE b\notin E\follows X_a\noncausal X_b\wrt{\ffilt{V}}$,
\item
$a\UE b\notin E\follows X_a\noncorr X_b\wrt{\ffilt{V}}$;
\end{romanlist}
\item[(LC)]
the {\em local Granger-causal Markov property}
with respect to $G$ if for all $a\in V$
\begin{romanlist}
\item
$X_{V\without(\parent{a}\cup\{a\})}\noncausal X_a\wrt{\ffilt{V}}$,
\item
$X_{V\without(\neighbour{a}\cup\{a\})}\noncorr X_a\wrt{\ffilt{V}}$;
\end{romanlist}
\item[(BC)]
the {\em block-recursive Granger-causal Markov property}
with respect to $G$ if for all subsets $A$ of $V$
\begin{romanlist}
\item
$X_{V\without(\parent{A}\cup A)}\noncausal X_A\wrt{\ffilt{V}}$,
\item
$X_{V\without(\neighbour{A}\cup A)}\noncorr X_A\wrt{\ffilt{V}}$.
\end{romanlist}
\end{list}
\end{definition}

Similarly, if $P$ is a probability kernel from $\rnum^{V\times\nnum}$ to
$\rnum^{V}$, we say that $P$ satisfies the pairwise, the local, or the
block-recursive Granger-causal Markov property with respect
to a graph $G$ whenever the same is true for every stationary process
$X_V$ with $X_V\sim P$.

\begin{figure}
\centerline{\includegraphics[width=\textwidth]{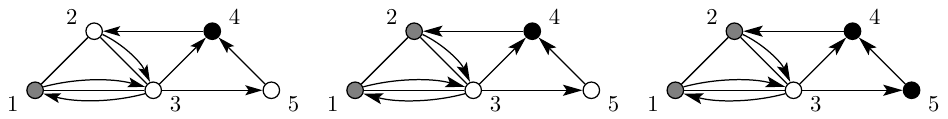}}
\caption{Encoding of relations $X_A\noncausal X_B\wrt{\ffilt{X}}$
by the (a) pairwise, (b) local, and (c) block-recursive
Granger-causal Markov property ($A$ and $B$ are indicated by grey and
black nodes, respectively).}
\label{fig-markovprop}
\end{figure}

\begin{example}
To illustrate the various Granger-causal Markov properties, we
consider the graph $G$ in Figure \ref{fig-markovprop}. Suppose that
a stationary process $X_V$ satisfies the pairwise Granger-causal Markov
property with respect to this graph $G$. Then the absence of the edge
$1\DE 4$ in $G$ implies that $X_1$ is Granger-noncausal for $X_4$ with
respect to $\ffilt{V}$. Next, in the case of the local Granger-causal
Markov property, we find that the bivariate subprocess $X_{\{1,2\}}$ is
Granger-noncausal for $X_4$ with respect to $\ffilt{V}$ since vertex 4
has parents 3 and 5. Similarly, if $X_V$ obeys the block-recursive
Granger-causal Markov property, the graph encodes that $X_{\{1,2\}}$ is
Granger-noncausal for $X_{\{4,5\}}$ with respect to $\ffilt{V}$ since
$\parent{4,5}=\{3\}$.
\end{example}

The block-recursive Granger-causal Markov property obviously implies the other
two Granger-causal Markov properties and, thus, is the strongest of the
three Markov properties; similarly, the pairwise Granger-causal Markov
property clearly is the weakest of the three properties. The question arises
whether and under which conditions the three Granger-causal Markov properties
are equivalent. In the case of random vectors $Y_V=(Y_v)_{v\in V}$ with values
in $\rnum^V$, the various levels of Markov properties for graphical
interaction models are equivalent if the distribution of $Y_V$ satisfies
\begin{equation}
\label{vector intersection}
Y_A\indep Y_B\given Y_{C\cup D}\wedge
Y_A\indep Y_C\given Y_{B\cup D}\follows
Y_A\indep Y_{B\cup C}\given Y_D
\end{equation}
for all disjoints subsets $A$, $B$, $C$, and $D$ of $V$ \citep{pearlpaz87}.
A necessary and sufficient condition for this intersection property is
that the information common to $Y_{B\cup D}$ and $Y_{C\cup D}$ equals the
information provided by $Y_D$. More precisely, let $(\Omega,\falg,\prob)$
be the underlying probability space and let $\yclass_S$ be the
sub-$\sigma$-algebra generated by $Y_S$, $S\subseteq V$. Furthermore,
we denote the $\sigma$-algebra generated by $\yclass_S$ and the $\prob$-null
sets in $\fclass$ by $\overline{\yclass_S}$. Then the above intersection
property holds if and only if
$\overline{\yclass_{C\cup D}}\cap\overline{\yclass_{B\cup D}}=
\overline{\yclass_D}$ \citep{dawid80,fmr90}; we say that $\yclass_{C\cup D}$
and $\yclass_{B\cup D}$ are measurable separable conditionally on
$\yclass_D$. For more details on measurable separability we refer to
Appendix \ref{appendix condind} and the references therein.

In order to ensure validity of the intersection property in the time series
case, we impose the following condition:
\begin{list}{{}}{\usecounter{alphcount}
\labelwidth2.3em
\leftmargin2.8em\labelsep0.5em\topsep0.25em plus 0.5ex
\itemsep0.25em plus 0.5ex\parsep0em}
\item[{\upshape(S)}]
for all subsets $A,B,C$ of $V$, $\ffilt{A}(t)$ and $\ffilt{B}(t)$ are
measurably separable conditionally on $\ffilt{A\cap B}(t)\vee\ffilt{C}(t-k)$
for all $k\in\nnum$ and $t\in\znum$.
\end{list}
Here, $\ffilt{A\cap B}(t)\vee\ffilt{C}(t-k)$ denotes the smallest
$\sigma$-algebra generated by $\ffilt{A\cap B}(t)\cup\ffilt{C}(t-k)$.
The condition implies that for every $\falg$-measurable random variable
$Z$ and all $t\in\znum$,
\begin{equation}
\label{intersection}
Z\indep\ffilt{A}(t)\given\ffilt{B\cup C}(t)\wedge
Z\indep\ffilt{B}(t)\given\ffilt{A\cup C}(t)\iff
Z\indep\ffilt{A\cup B}(t)\given\ffilt{C}(t).
\end{equation}

In the case of random vectors $Y_V$, a commonly used sufficient condition for
the intersection property and thus for conditional measurable separability is
that the joint distribution of $Y_V$ is absolutely continuous with respect to
some product measure and has a positive and continuous density
\citep[e.g.,][Prop.~3.1]{SL96}. The following result establishes a similar
condition in terms of conditional distributions for the time series case;
it requires an additional regularity condition on partial
tail-$\sigma$-algebras \citep{fm82,fmr90}.

\begin{proposition}
\label{intersectprop}
Let $X_V=\big(X_V(t)\big)_{t\in\znum}$ be a strictly stationary stochastic
process on some probability space $(\Omega,\falg,\prob)$ taking values in
$\rnum^{V}$ and suppose the following two conditions hold:
\begin{list}{{}}{\usecounter{alphcount}
\labelwidth2.3em
\leftmargin2.8em\labelsep0.5em\topsep0.25em plus 0.5ex
\itemsep0.25em plus 0.5ex\parsep0em}
\item[{\upshape(P)}]
the conditional distribution $\prob^{X_V(t+1)|\pastx_V(t)}$, $t\in\znum$,
has a regular version that is almost surely absolutely continuous with
respect to some product measure $\nu$ on $\rnum^{|V|}$ with
$\nu$-a.e.~positive and continuous density;
\item[{\upshape(M)}]
for all $A\subseteq V$ and $t\in\znum$
\[
\lcap_{k\in\nnum}\big(\overline{\ffilt{A}(t)}
\vee\overline{\ffilt{V\without A}(t-k)}\big)
=\overline{\ffilt{A}(t)}.
\]
\end{list}
Then the process $X_V$ satisfies condition (S).
\end{proposition}

For an interpretation of condition (M), we note that it is equivalent
to
\[
\lim_{k\to\infty}\mean\big(Z\given\ffilt{A}(t)\vee\ffilt{B}(t-k)\big)
=\mean\big(Z\given\ffilt{A}(t)\big).
\]
for all random variables $Z$ and subsets $A,B\subseteq V$
\citep{chamberlain82}. Thus condition (M) implies that the process $X_V$ is
conditionally weakly mixing. For many types of non-linear time series stronger
forms of mixing---but not conditional mixing---have been established
\citep[e.g.,][]{doukhan94,fanyao03}. We believe that the above condition
of conditional mixing is satisfied by most stationary time series models
but a discussion of this is beyond the scope of this paper.

The intersection property now allows us to derive the following relations
among the three Granger-causal Markov properties.

\begin{theorem}
\label{basicMP}
Suppose that $X_V$ satisfies condition (S). Then the three
Granger-causal Markov properties (BC), (LC), and (PC)
are related by the following implications:
\[
\mathrm{(BC)}\follows\mathrm{(LC)}\iff\mathrm{(PC)}.
\]
Furthermore, if $X_V$ additionally satisfies the composition property
\begin{equation}
\label{cond block}
X_A\noncausal X_B\wrt{\ffilt{V}}\iff
X_A\noncausal X_b\wrt{\ffilt{V}}\quad\forall\,b\in B,
\end{equation}
then the three Granger-causal Markov properties (BC), (LC), and (PC)
are equivalent.
\end{theorem}

The theorem shows that, similarly as in the case of chain graph models
with the Andersson-Madigan-Perlman (AMP) Markov property \citep{amp1}, the
pairwise and the local Granger-causal Markov property are in general
not sufficiently strong to encode all Granger-causal relationships that
hold among the components of a multivariate time series with respect to
full information $\ffilt{V}$. This suggests to specify graphical
time series models in terms of the block-recursive Granger-causal Markov
property.

\begin{definition}[Graphical time series model]
\label{graphmodel}
Let $G$ be a mixed graph and let $\pclass_G$ be a statistical time series
model given by a family of probability kernels $P\in\pclass_G$ from
$\rnum^{V\times\nnum}$ to $\rnum^V$. Then $\pclass_G$ is said to be
a {\em graphical time series model} associated with the graph $G$ if, for all
$P\in\pclass_G$, the distribution $P$ satisfies the block-recursive
Granger-causal Markov property with respect to $G$.
\end{definition}

The three Granger-causal Markov properties considered so far encode only
Granger noncausality relations with respect to the complete information 
$\ffilt{V}$. The discussion of phenomena such as spurious causality
\citep[e.g.,][]{hsiao82,eichlerbrain}, however,
requires also the consideration of Granger-causal relationships with respect
to partial information sets, that is, with respect to filtrations $\ffilt{S}$
for subsets $S$ of $V$. To this end, we introduce in Section \ref{sect:globalMP}
a global Granger-causal Markov property that more generally relates pathways in
a graph to Granger-causal relations among the variables, and we establish,
under condition (S), its equivalence to the block-recursive
Granger-causal Markov property; this shows that the block-recursive
Granger-causal Markov property is indeed sufficiently rich to describe
the dynamic dependence structure in  multivariate time series.

Before we continue our discussion of Markov properties in Section
\ref{sect:globalMP}, we illustrate the introduced concept of graphical
time series models by a few examples.

\section{Examples} \label{sect:examples}

In the previous section, graphical time series models have been defined
in terms of the block-recursive Granger-causal Markov property. For many
time series models, however, condition \eqref{cond block} in Theorem
\ref{basicMP} holds, and, hence, the pairwise, the local, and the
block-recursive Granger-causal Markov property are equivalent. This
enables us to derive the constraints on the parameters from the
pairwise or the local Granger-causal Markov property.

There are no simple conditions known that are both necessary and sufficient
for \eqref{cond block}. The following proposition lists some sufficient
conditions that cover many examples, as will be shown subsequently.

\begin{proposition}
\label{suffcond}
Suppose that $X_V$ satisfies condition (S) and one of the
following conditions:
\begin{romanlist}
\item
$X_V$ is a Gaussian process;
\item
$X_v(t+1)$, $v\in V$, are mutually contemporaneously independent,
that is, the joint conditional distribution factorizes as
\[
\prob^{X_V(t+1)\given\pastx_{V}(t)}
=\otimes_{v\in V}\prob^{X_v(t+1)\given\pastx_{V}(t)}\quad\forall t\in\znum;
\]
\item
$X_V(t+1)$ depends on its past only in its conditional mean, that is,
\[
X(t+1)-\mean\big[X(t+1)\given\ffilt{V}(t)\big]\indep\ffilt{V}(t)\quad
\forall t\in\znum.
\]
\end{romanlist}
Then the three Granger-causal Markov properties (BC), (LC), and (PC)
are equivalent.
\end{proposition}

We note that processes satisfying condition (ii) can be described by
directed graphs, that is, graphs without undirected edges. Thus the
proposition implies that for directed graphs the pairwise and the
block-recursive Granger-causal Markov property are always equivalent.

\subsection{Nonlinear autoregressive models}

As a first example, we consider the general class of multivariate
nonlinear autoregressive models given by
\[
X_V(t)=f_V\big(X_V(t-1),\ldots,X_V(t-p)\big)+\veps_V(t),
\]
where $f_V$ is an $\rnum^V$-valued Borel measurable function on
$\rnum^{p\times V}$ and $\veps_V=\big(\veps_V(t)\big)_{t\in\znum}$
is a sequence of independent and identically distributed zero mean
random vectors with density $q_V$ and such that $\veps(t)$ is independent
of $\ffilt{V}(t-1)$. Such models have been considered by many authors;
in particular, conditions on $f_V$ and $q_V$ that guarantee
geometric ergodicity and thus strong mixing of $X_V$ have been
established \citep[e.g.,][]{doukhan94,lujiang01,liebscher05}.
We note, however, that currently there are no conditions known
that ensure the conditional mixing condition (M). An exception are Gaussian
autoregressive processes that will be briefly discussed below.

For the general class of multivariate nonlinear autoregressive models, the
constraints imposed by a graph $G$ are
best formulated in terms of the local Granger-causal Markov
property. More precisely, $X_V$ satisfies the local Granger-causal
Markov property with respect to $G$ if for all $a\in V$
\begin{Alist}{L}
\item
$f_{a}\big(X_V(t-1),\ldots,X_V(t-p)\big)=
f_{a}\big(X_{\parent{a}\cup\{a\}}(t-1),\ldots,X_{\parent{a}\cup\{a\}}(t-p)\big)$;
\item
$q_V$ factorizes as
$q_{V}(z_V)=g_{a}\big(z_{\neighbour{v}\cup\{a\}}\big)\,
h_a(z_{V\without\{a\}})$.
\end{Alist}
The second condition implies $\veps_a(t)\indep
\veps_{V\without(\neighbour{a}\cup\{a\})}(t)\given\veps_{\neighbour{a}}$
which is equivalent to $X_a$ and $X_{V\without(\neighbour{a}\cup\{a\})}$
being contemporaneously conditionally independent with respect to $\ffilt{V}$
as required by the local Granger-causal Markov property.
Since $X_V(t)$ depends on its past $\pastx_V(t-1)$ only in its
conditional mean, it follows from Theorem \ref{basicMP} and Proposition
\ref{suffcond}(iii) that the local and the block-recursive Granger-causal
Markov properties are equivalent, that is, the above conditions on
$f_V$ and $q_V$ define indeed a graphical nonlinear autoregresssive model
of order $p$ associated with the graph $G$.

The general class of multivariate nonlinear autoregressive models
covers many interesting and important models, of which we discuss only
the following three.
\begin{alphlist}
\item
\label{VAR example}
{\em Vector autoregressive (VAR) model:}
Suppose that $X_V$ is a stationary Gaussian process given by
\begin{equation}
\label{VAR-def}
X_V(t)=\lsum_{u=1}^p\Phi(u)\,X_V(t-u)+\veps(t),
\qquad\veps(t)\iid\normal(0,\Sigma),
\end{equation}
where $\Phi(u)$ are $V\times V$ matrices and the variance matrix $\Sigma$
is non-singular with inverse $K=\Sigma^{-1}$. Then $X_V$ satisfies the
pairwise Granger-causal Markov property with respect to a graph $G=(V,E)$ if
for all $a\neq b$
\begin{romanlist}
\item
$a\DE b\notin E\follows \Phi_{ba}(u)=0\quad\forall u=1,\ldots,p$;
\item
$a\UE b\notin E\follows K_{ab}=K_{ba}=0$.
\end{romanlist}
Thus, the graphical VAR model of order $p$ associated
with the graph $G$, denoted by VAR($p$,$G$), is given by all
stationary VAR($p$) processes whose parameters are constrained to zero
according to the conditions (i) and (ii).

Furthermore, let
$f(\lam)=(2\pi)^{-1}\,\Phi(e^{-\im\lam})^{-1}\,\Sigma\,\Phi(e^{-\im\lam})'{}^{-1}$,
$\lam\in[-\pi,\pi]$, be the spectral density matrix of $X_V$, where
$\Phi(z)=I_V-\Phi(1)\,z-\ldots-\Phi(p)\,z^p$ and $I_V$ is the
$V\times V$ identity matrix. Then, if the eigenvalues of $f(\lam)$ are
bounded and bounded away from zero uniformly for all $\lam\in[-\pi,\pi]$,
the process $X_V$ satisfies the separability condition (S)
\citep[Lemma A.2]{eichlerpathdiagr}.
\item
{\em Self-exciting threshold autoregressive (SETAR) model:}
A stochastic process $X_V$ is said to follow a multivariate SETAR model
\citep[e.g.,][]{tong93,arnold01} if for each $a\in V$
\[
X_a(t)=\lsum_{u=1}^{p}\lsum_{b\in V}
\phi^{(n)}_{ab}(u)\,X_b(t-u)+\veps_a(t)
\qquad\text{if }X_a(t-d)\in I_{a,n},
\]
where $\{I_{a,1},\ldots,I_{a,N}\}$ is a partition of $\rnum$, and
$\veps_V(t)\iid Q_V$, say. Then $X_V$ obeys the local Granger-causal
Markov property with respect to a graph $G=(V,E)$ if, for $a\neq b$,
$\phi^{(n)}_{ab}(u)=0$ for all $n=1,\ldots,N$ and $u=1,\ldots,p$
whenever $b\DE a\notin E$ and $Q_V$ has density $q_V$ satisfying
condition (L2).
\item
{\em Nonparametric additive autoregressive model}:
A very useful class of nonparametric autoregressive models, which avoid the
``curse of dimensionality'', are the additive models given by
\[
X_a(t)=\lsum_{u=1}^{p}\lsum_{b\in V}f_{ab}^{(u)}\big(X_b(t-u)\big)
+\veps_a(t),\quad a\in V,\,t\in\znum,
\]
where $f^{(u)}_{ab}$ are real-valued functions on $\rnum$. Here,
condition (L1) obviously is equivalent to that the functions $f_{ab}^{(u)}$,
$u=1,\ldots,p$, are constant whenever $a\neq b$ and the edge $b\DE a$ is
missing in the graph $G$.
\end{alphlist}

\subsection{Multivariate ARCH processes}
\label{section:arch}

Another important class of nonlinear time series models are
the autoregressive conditional heteroscedasticity (ARCH) model
and its various subsidiaries, which have been developed for
modelling the time-varying volatility exhibited by many
financial time series. A stationary stochastic process $X_V$ is
said to follow a multivariate ARCH($q$) process if its conditional
mean $\mean(X(t)\given\ffilt{X}(t-1))$ is zero and the conditional
covariance matrix is of the form
\[
\mean\big(X(t)X(t)'\given\ffilt{V}(t-1)\big)=\Sigma(t)=
g_{VV}\big(X_V(t-1),\ldots,X_V(t-p)\big).
\]
For an overview of multivariate ARCH models
we refer to \citet{bollerslevarch} and \citet{archmodelsbook};
sufficient conditions ensuring existence and strong mixing of
such processes can be found, for instance, in
\citet{lujiang01}, \citet{carrasco02}, and \citet{liebscher05}.

One key issue in the specification of multivariate ARCH models
is the restriction of the number of parameters involved, which
in a general setting can be very large. Various parametrisations
that allow different levels of complexity have been suggested.
Here the graphical modelling approach can help to achieve
a further reduction of the number of parameters.

In the following, we consider stochastic processes $X_V$
with conditional distribution $\normal(0,\Sigma(t))$ and
formulate the constraints defining a graphical ARCH($q$) model
associated with a graph $G=(V,E)$ for three different
parametrisations of $\Sigma(t)$.
\begin{romanlist}
\item
{\itshape Constant conditional correlations:}
The constant conditional correlation model of \citet{bollerslev90}
provides the most parsimonious parametrisation of $\Sigma(t)$.
The conditional variances are given by
\begin{align*}
\sigma_{aa}(t)
&=\sigma_{aa}^0+\lsum_{u=1}^q\lsum_{k\in\parent{a}\cup\{a\}}\alpha^a_k(u) X_k(t-u)^2,
\end{align*}
whereas the conditional covariances are
determined by the set of equations
\begin{align*}
\sigma_{ab}(t)
&=\sigma_{aa}(t)^{1/2}\sigma_{bb}(t)^{1/2}\rho_{ab}&&\text{ if }a\UE b\in E,\\
K_{ab}(t)&=0&&\text{ if }a\UE b\notin E.
\end{align*}
Here $K(t)=\Sigma(t)^{-1}$ is the inverse conditional covariance matrix.
\item
{\itshape Constant conditional correlations with interaction:}
In this parametrisation the conditional variance $\sigma_{aa}(t)$
additionally depends on interaction terms $X_k(t-u)X_l(t-u)$ 
if $k$ and $l$ are
both parents of $a$. Thus the conditional variance can be written as
\begin{align*}
\sigma_{aa}(t)
&=\sigma_{aa}^0+\lsum_{u=1}^q\lsum_{\substack{k,l\in\parent{a}\cup\{a\}:k<l}}
\alpha^a_{kl}(u)X_k(t-u)X_l(t-u).
\end{align*}
The entries $\sigma_{ab}(t)$ have the same form as in (i).
\item
{\itshape Vector ARCH model:}
In the general vector ARCH model due to \citet{kraftengle82}, also the correlation between the components of $X(t)$ may depend on the past values of $X$.
This leads to conditional covariances $\sigma_{ab}(t)$, $a\leq b$,
of the form 
\begin{align*}
\sigma_{ab}(t)
&=\sigma_{ab}^0+\lsum_{u=1}^q
\lsum_{\substack{k,l\in P_{ab}:k<l}}
\alpha^{ab}_{kl}(u)X_k(t-u)X_l(t-u)
\end{align*}
if $a=b$ or $a\UE b\in E$, where
$P_{ab}=(\parent{a}\cup\{a\})\cap(\parent{b}\cup\{b\})$,
while the conditions $K_{ab}(t)=0$
for $a\neq b$ and $a\UE b\notin E$ remain unchanged.
\end{romanlist}

For the constant conditional correlation models it is easy to
derive conditions to ensure that the conditional covariances are
positive definite almost surely for all $t$. In contrast, such
conditions are difficult to impose and verify for the vector
ARCH model. Therefore \citet{englekroner} suggested
an alternative representation for the multivariate ARCH($q$) model
in which $\Sigma(t)$ is guaranteed to be positive definite
almost surely for all $t$. In this so-called BEKK
representation\footnote{This is named after Baba, Engle, Kraft and Kroner,
the authors of an earlier version of the paper \citep[cf][]{baba90}.},
the conditional covariances of a graphical ARCH model are
parametrised by
\[
\sigma_{ab}(t)=
\sigma^{0}_{ab}
+\lsum_{n=1}^N\lsum_{u=1}^q
\lsum_{\substack{k,l\in P_{ab}:k<l}}
\alpha^{(n)}_{ka}(u)\alpha^{(n)}_{lb}(u)X_k(t-u)X_l(t-u).
\]
In this form it is immediately clear that if $\sigma_{ab}(t)$
depends on the past of $X_k$ then at least one of the conditional
variances $\sigma_{aa}(t)$ and $\sigma_{bb}(t)$ must also depend
on $X_k$. Although less obvious the same can be shown for the vector
ARCH model in the original parametrisation noting that the conditional
covariance matrix $\Sigma(t)$ must be positive definite. Hence
graphical vector ARCH models fulfill condition \eqref{cond block}.
For the constant conditional correlation model condition
\eqref{cond block} is trivially fulfilled.

\subsection{A binary time series model}
\label{binarymodel}

As an example with categorical data, we consider a binary time series
model that has been used for the identification of neural interactions
from neural spike train data \citep{brillinger88:biolcyb,brillinger88:biomedeng}.
Suppose that the data consist of the recorded spike trains for a set of
neurons, that is, of the sequences of firing times $(\tau_{v,n})_{n\in\nnum}$
for neurons $v\in V$, and let $X_v$ be the binary time series obtained by
setting $X_v(t)=1$ if neuron $v$ has fired in the interval $[t,t+1)$
and $X_v(t)=0$ otherwise. We assume that the hypothesized neural pathways between the observed neurons can be depicted by a purely directed graph $G$; in particular, we thus exclude the possibility that the dependences among the observed neurons are affected by unmeasured confounders. Then the interactions between the neurons can be modelled by the conditional probabilities
\begin{equation}
\label{neuralmodel1}
\prob\big(X_b(t)=1\big|\ffilt{V}(t-1)\big)
=\Phi\Big(\lsum_{a\in\parent{b}}U_{ba}(t)-\theta\Big),
\end{equation}
where $\Phi(x)$ denotes the normal cumulative function,
\begin{equation}
\label{neuralmodel2}
U_{ba}(t)=\lsum_{u=1}^{\gamma_b(t)}g_{ba}(u)\,X_a(t-u)
\end{equation}
measures the influence of process $a$ on process $b$, and
\[
\gamma_b(t)=\min\big\{u\in\nnum\big|X_b(t-u)=1\big\}
\]
is the time elapsed since the last event of process $X_b$. Furthermore,
we assume that the time unit has been chosen small enough such that there
are no interactions among the neurons within one time interval, and that,
consequently, the joint conditional probability factorizes as
\[
\prob\big(X_V(t)=x_V\big|\ffilt{V}(t-1)\big)
=\lprod_{v\in V}\prob\big(X_v(t)=x_v\big|\ffilt{V}(t-1)\big)
\]
for all $x_V\in\{0,1\}^{V}$. Then the pairwise and the block-recursive
Granger-causal Markov property are equivalent by Proposition \ref{suffcond}(ii)
and, thus, we can use the former for modelling dependences between the
processes. From \eqref{neuralmodel1} and \eqref{neuralmodel2}, it follows
that $X_a$ is Granger-noncausal for $X_b$ if and only if $g_{ba}(u)=0$ for
all $u\in\nnum$.

\subsection{Two counter examples}

Although condition \eqref{cond block} is satisfied by a wide variety
of time series models it does not hold generally. As an example, we
consider a simple nonlinear ARCH model $X_V$ with conditional distributions
$X_V(t)|\ffilt{V}(t-1)\sim\normal\big(0,\Sigma(t)\big)$,
where the conditional covariance matrix $\Sigma(t)$ is given by
\begin{equation}
\label{ex-nonequiv}
\Sigma(t)=\begin{pmatrix}
1&\rho(t)&0\\
\rho(t)&1&0\\
0&0&1
\end{pmatrix}\quad\text{with}\quad
\rho(t)=\begin{mycases}
\rho&\text{ if }|X_3(t-1)|>c\\
0&\text{ otherwise}
\end{mycases}
\end{equation}
for some constants $\rho$ with $0<|\rho|<1$ and $c>0$.
Models of this type can be seen as a multivariate generalisation of the
qualitative threshold ARCH($1$) model of \citet{gourierouxmonfort}.

\begin{figure}
\centerline{\includegraphics[width=7cm]{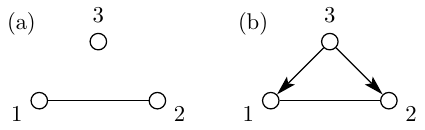}}
\caption{Illustration of non-equivalence of pairwise and block-recursive
Granger-causal Markov properties: the process with conditional variance
\eqref{ex-nonequiv}
satisfies the pairwise Granger-causal Markov property with respect
to the graphs in (a) and (b) whereas it satisfies the block-recursive
Granger-causal Markov property only with respect to the graph in (b).}
\label{fig-nonequiv}
\end{figure}

From the conditional covariance matrix, we find that, on the one hand,
the marginal conditional distributions of $X_v(t)$ given $\ffilt{V}(t-1)$ are
standard normal and, thus, do not depend on $\ffilt{V}(t-1)$. This implies
that the process $X_V$ satisfies the pairwise Granger-causal Markov property
with respect to the graph (a) in Figure \ref{fig-nonequiv}. On the other hand,
$X_k$ Granger-causes the subprocess $(X_1,X_2)$ since the bivariate conditional
distribution of $\big(X_1(t),X_2(t)\big)$ depends on the value of $X_3(t-1)$
through the conditional correlation $\rho(t)$. Thus $X_V$ obeys the
block-recursive Granger-causal Markov property with respect to the graph (b)
in Figure \ref{fig-nonequiv}, but not with respect to the graph (a).

We note that the example can be easily generalized by considering
models where the conditional variances $\var\big(X_a(t)\big|\ffilt{V}(t-1)\big)$,
$a\in V$, and the conditional correlation matrix
$\corr\big(X_V(t),X_V(t)\big|\ffilt{V}(t-1)\big)$ are modelled separately
as functions of the past values $X_V(t-1),\ldots,X_V(t-p)$.
\smallskip

Next, consider the trivariate process $X_V$ given by
\[
X_1(t)=f\big(X_2(t-1)\big)+\veps(t),\quad
X_2(t)=g\big(X_3(t-1)\big),\quad
X_3(t)=\eta(t),
\]
where $\veps(t)$ and $\eta(t)$ are independent sequences of
{\em i.i.d.~}random variables. Since
\[
\ffilt{1}(t)\vee\ffilt{2}(t)=\ffilt{1}(t)\vee\ffilt{3}(t-1),
\]
condition (S) is violated. Indeed, we find that neither $X_2$ nor $X_3$ Granger-cause $X_1$ with respect to the full filtration $\ffilt{V}$ whereas the bivariate process $(X_2,X_3)'$ is Granger-causal for $X_1$. Therefore, the pairwise and the local Granger-causal Markov property are not equivalent for this process.

\section{Global Markov properties}
\label{sect:globalMP}

The interpretation of graphs describing the dependence structure
of graphical models in general is enhanced by global Markov properties
that merge the notion of conditional independence with a purely graph
theoretical concept of separation allowing one to state whether two
subsets of vertices are separated by a third subset of vertices.
In this section, we show that the concept of $p$-separation
introduced by \citet{amp2} for chain graph models with the AMP Markov
property \citep{amp1} can be used to obtain global Markov properties in
the present context of graphical time series models.
Throughout this section we assume that condition (S) in Section
\ref{sect:graphmodel} holds.

\subsection{The global AMP Markov property}

We start with some further graphical terminology. Let $G=(V,E)$ be a
mixed graph. Then a {\em path} $\pi$ between two vertices $a$ and $b$
in $G$ is a sequence $\pi=\path{e_1,\ldots,e_n}$ of edges $e_i\in E$
such that $e_i$ is an edge between $v_{i-1}$ and $v_i$ for some sequence
of vertices $v_0=a,v_1,\ldots,v_{n}=b$. The vertices $a$ and $b$ are
the {\em end-points} of the path, while $v_1,\ldots,v_{n-1}$ are the
{\em intermediate points} on the path. Like \citet{koster02} we do
not require that the points $v_j$ on a path $\pi$ are distinct; this
means that paths in general may be self-intersecting.
A path $\pi$ in $G$ is called a {\em directed path} if it is of the
form $a\DE\ldots\DE b$ or $a\LDE\ldots\LDE b$. Similarly, if $\pi$
consists only of undirected edges it is called an {\em undirected path}.
Furthermore, a path $\tilde\pi$ is a subpath of $\pi$ if
$\tilde\pi=\path{e_i,e_{i+1},\ldots,e_{j-1},e_j}$ for some $1\leq i\leq j\leq n$.

An intermediate point $c$ on a path $\pi$ is said to be a
{\em $p$-collider} on the path if the edges preceding and suceeding
$c$ on the path either have both an arrowhead at $c$ or one has an arrowhead
at $c$ and the other is a line, i.e.~$\DE c \LDE$, $\DE c \UE$,
$\UE c \LDE$; otherwise the point $c$ is said to be a
{\em $p$-noncollider} on the path. Notice that this classification only
applies to the intermediate points of a path $\pi$; the end-points are
neither $p$-colliders nor $p$-noncolliders. We also note that a vertex can
take different roles in different positions on a path: for example, on the
path $1\DE 3\LDE 2\DE 3\DE 4$ in Figure \ref{fig-markovprop}, vertex $3$ appears both as an $p$-collider and an $p$-noncollider.

A path $\pi$ between vertices $a$ and $b$ is said to be {\em $p$-connecting}
given a set $S$ if
\begin{romanlist}
\item
every $p$-noncollider on the path is not in $S$, and
\item
every $p$-collider on the path is in $S$,
\end{romanlist}
otherwise we say the path is {\em $p$-blocked} given $S$. In graphs
encoding Markov properties of variables, $p$-connecting paths
are exactly those paths inducing associations between the variables;
conversely, if there are no $p$-connecting paths the corresponding
variables are independent. This leads to the following definition.

\begin{definition}[\boldmath$p$-separation]
Two vertices $a$ and $b$ in a mixed graph $G$ are $p$-separated
given a set $S$ if all paths between $a$ and $b$ are $p$-blocked
given $S$. Similarly, two sets $A$ and $B$ in $G$ are said to be
$p$-separated given $S$ if, for every pair $a\in A$ and  $b\in B$,
$a$ and $b$ are $p$-separated given $S$. This will be denoted by
$A\sepp B\given S$.
\end{definition}

We note that the above conditions for $p$-separation are simpler
than those in \citet{amp2} due to the fact that we consider
the larger class of all possibly self-intersecting paths.
The equivalence of the two notions of $p$-separation
is shown in Appendix \ref{psepequivalence}.
The following results show that the concept of $p$-separation
can be applied to graphs encoding dynamic relationships in multivariate
time series and allows reading off conditional independences among the
stochastic processes that are represented by the vertices in the graph.

\begin{lemma}
\label{globalMPpast}
Suppose that $X_V$ satisfies the block-recursive Granger-causal Markov property
with respect to the graph $G$. Then, for any disjoint subsets $A$, $B$, and $S$
of $V$, we have
\[
A\sepp B\given S\follows
\ffilt{A}(t)\indep\ffilt{B}(t)\given\ffilt{S}(t)\qquad\forall t\in\znum.
\]
\end{lemma}

Letting $t$ tend to infinity, we can translate $p$-separation in the graph
into conditional independence statements for complete subprocesses.
For this, we define $\ffilt{S}(\infty)=\vee_{t\in\znum}\ffilt{S}(t)$
as the $\sigma$-algebra generated by the subprocess $X_S$.

\begin{theorem}
\label{global AMP}
Suppose $X_V$ satisfies the block-recursive Granger-causal Markov property
with respect to the graph $G$. Then, for any disjoint subsets $A$, $B$,
and $S$ of $V$, we have
\[
A\sepp B\given S\follows
\ffilt{A}(\infty)\indep\ffilt{B}(\infty)\given\ffilt{S}(\infty).
\]
We say that $X$ satisfies the {\em global AMP Markov property (GA)} with
respect to $G$.
\end{theorem}

\begin{figure}
\centerline{\includegraphics[width=\textwidth]{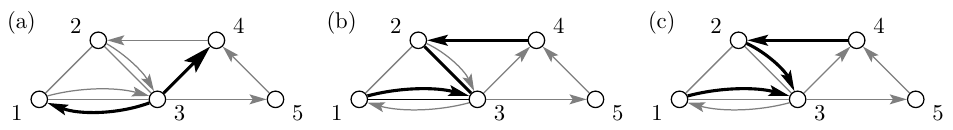}}
\caption{Illustration of global AMP Markov property (paths are
marked by bold lines):
(a) path between $1$ and $4$ that is $p$-connecting given $S\subseteq\{2,5\}$;
(b) path between $1$ and $4$ that is $p$-connecting given $S=\{2,3\}$
(or $\{2,3,5\}$);
(c) path between $1$ and $4$ that is $p$-connecting given $S=\{3,5\}$
(or $\{3\}$).}
\label{fig-ampmarkovprop}
\end{figure}

\begin{example}
\label{ex-ampmarkov}
For an illustration of the global AMP Markov property, we consider again
the graph $G$ in Figure \ref{fig-markovprop}. In this graph, vertices
$1$ and $4$ are not adjacent. Nevertheless, it can be shown that the
two vertices cannot be $p$-separated by any set $S\subseteq\{2,3,5\}$:
firstly, the path $1\LDE 3\DE 4$ is
$p$-connecting given a set $S$ unless the set $S$ contains the vertex $3$
(Fig.~\ref{fig-ampmarkovprop} a). Secondly, the path $1\DE 3\UE 2\LDE 4$
is $p$-connecting given $S$ whenever both intermediate points $2$ and
$3$ belong to $S$ (Fig.~\ref{fig-ampmarkovprop} b). Finally, the path
$1\DE 3\LDE 2\LDE 4$ is $p$-connecting given $S$ if $S$ contains vertex
$3$ but not $2$ (Fig.~\ref{fig-ampmarkovprop} c).
Thus, if $X_V$ is a stationary process that obeys the block-recursive
Granger-causal Markov property with respect to $G$, then the graph $G$
does not encode that $X_1$ and $X_4$ are conditionally independent given
$X_S$ regardless of the choice of $S\subseteq\{2,3,5\}$.

Similarly, it can be shown that vertices $1$ and $5$ are $p$-separated
given $S=\{3,4\}$: every path between $1$ and $5$ that contains the
edge $3\DE 5$ or the subpath $3\DE 4\LDE 5$ is $p$-blocked by vertex $3$.
All other paths between $1$ and $5$ contain the subpath $2\LDE 4\LDE 5$
and, thus, are blocked by vertex $4$. It follows that for every process
$X_V$ that satisfies the block-recursive Granger-causal Markov property
with respect to $G$ the components $X_1$ and $X_5$ are conditionally
independent given $X_{\{3,4\}}$.
\end{example}

\subsection{The global Granger-causal Markov property}

In this section, we apply the concept of pathwise separation to
the problem of deriving general Granger noncausality relations
from mixed graphs. To motivate the approach, we firstly consider
the graphical VAR(1) model of all trivariate stationary processes
$X_V=(X_1,X_2,X_3)$ given by
\begin{equation}
\label{trivariate-formula}
\begin{split}
X_1(t)&=\phi_{11}\,X_1(t-1)+\phi_{12}\,X_2(t-1)+\veps_1(t),\\
X_2(t)&=\phi_{22}\,X_2(t-1)+\phi_{23}\,X_3(t-1)+\veps_2(t),\\
X_3(t)&=\phi_{33}\,X_3(t-1)+\veps_3(t)
\end{split}
\end{equation}
for $t\in\znum$ with independent and standard normally distributed
errors $\veps_V(t)$, $t\in\znum$. The associated graph $G$ that encodes the
restrictions imposed on the model consists simply of the path
$3\DE 2\DE 1$, which is $p$-connecting given the empty set.
This indicates that the components $X_1$ and $X_3$ are, in general, not
independent in a bivariate analysis. However, an intuitive interpretation
of the directed path $3\DE 2\DE 1$ suggests that $X_3$ Granger-causes $X_1$
but not vice versa if only the bivariate process $X_{\{1,3\}}$ is considered.
Indeed, the block-recursive Granger-causal Markov property implies that
$X_3(t+1)\indep\ffilt{\{1,2\}}(t)\given\ffilt{\{3\}}(t)$,
from which it follows by decomposition (see Appendix A) that
$X_1$ is Granger-noncausal for $X_3$ with respect to $\ffilt{\{1,3\}}$.
Obviously, the $p$-separation criterion is too strong for establishing
this Granger-noncausality relationship between $X_3$ and $X_1$ since it
requires that all paths between the two vertices are $p$-blocked
whereas it seems sufficient that only certain paths, namely those
ending with an arrowhead at vertex $3$, are $p$-blocked.

This suggests the following definitions.
A path $\pi$ between two vertices $a$ and $b$ in $G$ is said to
be {\em $b$-pointing}\footnote{In the literature, a path with this property is also termed a path into $b$.} if it has an arrowhead at the endpoint
$b$. More generally, a path $\pi$ between two
disjoint subsets $A$ and $B$ is said to be {\em $B$-pointing} if
it is $b$-pointing for some $b\in B$.

For the derivation of contemporaneous conditional independences,
we also need to consider paths with arrowheads at both endpoints;
such paths $\pi$ will be called {\em bi-pointing}. Furthermore, let
$\pi=\path{\pi_1,\ldots,\pi_n}$ be a composition of paths $\pi_i$ that
are undirected or bi-pointing. Then $\pi$ is said to be an
{\em extended bi-pointing} path. In particular, this implies that
any undirected or bi-pointing path is also an extended bi-pointing path;
similarly, the composition $\pi=\path{\pi_1,\pi_2}$ of two extended
bi-pointing paths $\pi_i$ is again extended bi-pointing.
Moreover, every extended bi-pointing path $\pi$ is of the
form $\pi=\path{u_1,\beta,u_2}$ for some paths $u_1$,
$u_2$, and $\beta$ of possibly length zero, where $u_1$ and $u_2$ are
undirected paths and $\beta$ is a bi-pointing path (hence the term
`extended bi-pointing'). 
With these definitions, we define the following global Granger-causal
Markov property, which gives a path-oriented criterion for deriving
general Granger noncausality relations from a mixed graph.

\begin{definition}[Global Granger-causal Markov property]
Let $X_V$ be a stationary process and let $G=(V,E)$ be a mixed graph. Then
$X_V$ satisfies the {\em global Granger-causal Markov property} (GC) with
respect to $G$ if, for all disjoint subsets $A$, $B$, and $S$ of $V$, the
following conditions hold:
\begin{romanlist}
\item
if every $B$-pointing path in $G$ between $A$ and $B$ is $p$-blocked
given $S\cup B$ then
\[
X_A\noncausal X_B\wrt{\ffilt{A\cup B\cup S}};
\]
\item
if every extended bi-pointing path in $G$ between $A$ and $B$ is $p$-blocked
given $A\cup B\cup S$ then
\[
X_A\noncorr X_B\wrt{\ffilt{A\cup B\cup S}}.
\]
\end{romanlist}
\end{definition}

From the definition, it is immediately clear by setting $S=V\without(A\cup B)$
that the global Granger-causal Markov property entails the block-recursive
Granger-causal Markov property. The following theorem shows that in fact,
under condition (S), the two Granger-causal Markov properties are
equivalent; thus, the global Granger-causal Markov property may be employed
to discuss the dynamic relationships implied by a graphical time series model
defined in terms of the block-recursive Granger-causal Markov property.

\begin{theorem}
\label{global-markov}
Let $X_V$ be a stationary process and let $G=(V,E)$ be a mixed graph. Then
$X_V$ satisfies the block-recursive Granger-causal Markov property
with respect to $G$ if and only if $X_V$ satisfies the
global Granger-causal Markov property with respect to $G$.
\end{theorem}

As a consequence of the global Granger-causal Markov property,
we find that $p$-separation in the graph implies Granger noncausality
in both directions and contemporaneous conditional independence
of the variables.

\begin{corollary}
\label{psepgranger}
Suppose that the
process $X_V$ satisfies the block-recursive Granger-causal Markov property
with respect to a mixed graph $G$. For disjoint subsets $A$, $B$, and $S$
of $V$, if $A$ and $B$ are $p$-separated given $S$, then
\[
X_A\noncausal X_B\wrt{\ffilt{A\cup B\cup S}},\quad
X_B\noncausal X_A\wrt{\ffilt{A\cup B\cup S}},\text{ and }\,
X_A\noncorr X_B\wrt{\ffilt{A\cup B\cup S}}.
\]
\end{corollary} 

The following corollary summarizes the relationships between the various
Markov properties for graphical time series models.

\begin{corollary}
The various Granger-causal Markov properties are related as follows:
\begin{align*}
\mathrm{(GC)}\iff&\mathrm{(BC)}\follows\mathrm{(LC)}\iff\mathrm{(PC)}.
\end{align*}
Furthermore, we have $\mathrm{(BC)}\follows\mathrm{(GA)}$.
If additionally condition \eqref{cond block} holds, then the four
Granger-causal Markov properties (PC), (LC), (BC), and (GC) are
equivalent.
\end{corollary}

\begin{proof}
The corollary summarizes Theorems \ref{basicMP}, \ref{global AMP}, and \ref{global-markov}.
\end{proof}

\begin{example}
\label{expsep}
For an illustration, we again consider a stationary time series $X_V$
satisfying the block-recursive Granger-causal Markov property with respect
to the graph $G$ in Figure \ref{fig-markovprop}. In Example \ref{ex-ampmarkov},
we have seen that vertices $1$ and $4$ are not $p$-separated given $S=\{3\}$,
that is, $X_1$ and $X_4$ are in general not conditionally independent given
$X_3$. We now employ the global Granger-causal Markov property to examine
the dynamic relationships between the components $X_1$ and $X_4$ further.

We start by examining the $4$-pointing paths between $1$ and $4$. Straightforward considerations show that all $4$-pointing paths end with either $3\DE 4$, $3\DE 5\DE 4$, or $2\LDE 4\LDE 5\DE 4$; three instances of such paths are depicted in Figure \ref{fig-psep1}. The paths ending with either $3\DE 4$ or $3\DE 5\DE 4$ are clearly $p$-blocked by vertex $3$ whereas the paths ending with $2\LDE 4\LDE 5\DE 4$ are $p$-blocked by vertex $4$.
It follows that every $4$-pointing paths between $1$ and $4$ is $p$-blocked by $\{3,4\}$ and thus $X_1$ does not Granger-cause $X_4$ with
respect to $\ffilt{\{1,3,4\}}$.

\begin{figure}
\centerline{\includegraphics[width=\textwidth]{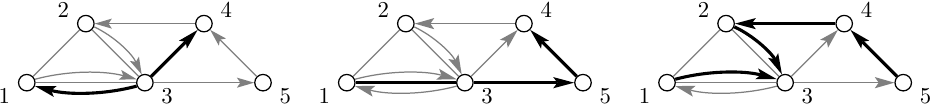}}
\caption{Illustration of global Granger-causal Markov property:
Three $4$-pointing paths (solid lines) between $1$ and
$4$ that are $p$-blocked by the set $\{3,4\}$.}
\label{fig-psep1}
\end{figure}

Similarly, we can examine all extended bi-pointing paths between
vertices $1$ and $4$ to show that $X_1$ and $X_4$ are contemporaneously
conditionally independent with respect to $\ffilt{\{1,3,4\}}$.
Figure \ref{fig-psep2} shows three examples of such paths: the first
two are $p$-blocked by vertex $3$ (notice that on the second path, the vertex
$3$ is once a $p$-collider and once a $p$-noncollider) whereas the
last path is $p$-blocked by vertices $3$ and $4$. For similar reasons
as above, these three paths are exemplary for all extended bi-pointing
paths between $1$ and $4$, and we conclude that $X_1$ and $X_4$ are
indeed contemporaneously conditionally independent with respect to
$\ffilt{\{1,3,4\}}$.

\begin{figure}
\centerline{\includegraphics[width=\textwidth]{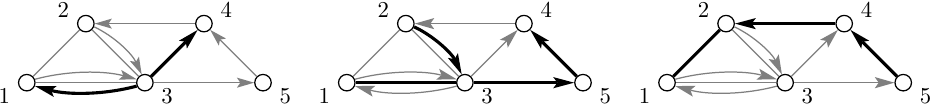}}
\caption{Illustration of global Granger-causal Markov property:
Three extended bi-pointing paths (solid lines) between $1$ and
$4$ that are $p$-blocked by the set $\{3,4\}$.}
\label{fig-psep2}
\end{figure}

Finally, we note that every $1$-pointing path between $4$ and $1$
must end with the directed edge $3\DE 1$. Since this edge has a
tail at vertex $3$, every such path must be $p$-blocked given
$S=\{1,3\}$, which implies that $X_4$ does not Granger-cause $X_1$
with respect to $\ffilt{\{1,3,4\}}$.
\end{example}

\section{Discussion}
\label{sect:discussion}
In this paper, we discussed a graphical modelling approach
for multivariate time series that is based on mixed graphs in
which each vertex  represents one complete component series
while the edges in the graph reflect possible dynamic interdependencies
among the variables of the process. The constraints imposed by the graphs
are formulated in terms of strong Granger noncausality and, thus, allow
modelling arbitrary non-linear dependencies. The graphical modelling
approach can help to reduce the number of parameters involved in modelling
high-dimensional non-linear time series while encoding the constraints
on the parameters in a simple graph, which is easy to visualize and allows
an intuitive understanding of the dependencies in the model.

We have shown that the interpretation of these graphs, which for
many models are built only from pairwise Granger noncausality relations,
is enhanced by so-called global Markov properties, which relate
separation properties of the graph to conditional independence or
Granger noncausality statements about the process. In this paper,
we have used the path-oriented concept of $p$-separation, which allows us
to attribute Granger-causal relationships among the variables
to certain pathways in the graphs.

Our objective has been to provide a general framework for modelling
the dynamic interdependencies in multivariate time series; in particular,
we focused on a simple graphical representation, which has been achieved
by representing each component of a multivariate time series by a single
vertex in the associated graph. The approach presented here, however, is not
the only possible, and since the first papers on the application of
graphical models in time series analysis \citep{lynggaardwalther,brill96},
there has been an increasing interest in the topic
\citep{stanghellini99,RD00,reale01,hsss,oxleyreale04,moneta05,
eichlerpathdiagr,eichlerhandbook}. All these approaches are basically
restricted to the analysis of linear interdependencies, and most of them
represent each variable at each time point by a separate vertex in the
associated graph. In the following, we briefly compare our approach
with alternative graphical representations and point out possible extensions.

\subsection*{\em Modelling processes of variables at separate time points}

A more detailed modelling of dependencies among the components
of a vector time series can be achieved by representing each random
variable $X_v(t)$ by a different vertex $v_t$, say, in a graph $G$.
This alternative approach has been discussed, for example, by
\citet{reale01}, \citet{hsss}, and \citet{moneta05}. On the one hand,
it leads to a more flexible class of graphical models and has the
advantage that many of the concepts and methods that have been
developed for the multivariate case carry over to the time series
case. On the other hand, the increased flexibility leads to
(sometimes much) larger graphs, which easily can become unwieldy
and difficult to interpret, and it clearly also aggrevates the
model selection problem. Moreover, the underlying graph for
such graphical time series models theoretically has infinitely
many vertices, and it is not immediately clear how to prune this
graph to a finite representation while preserving the Markov
properties. In contrast, Lemma \ref{finitepaths} provides a simple
local criterion that restricts the search for $p$-connecting paths
in the type of graphs considered in this paper.

Apart from these theoretical and practical issues, we think that
a high level of detail as provided by these models is not always
wanted nor always appropriate. We give two examples.
Firstly, \citet{baccala01}
proposed a frequency-domain approach for the discussion of
Granger-causal relationships based on the concept of partial
directed coherence. Although this approach still requires the
fitting of VAR models, the identification of interactions
is performed in the frequency-domain and hence only relations
on the level of Granger noncausality can be identified. The
results in \citet{baccala01} were summarized by path diagrams
associated with the identified VAR model as discussed in
\citet{eichlerpathdiagr}. Our approach of representing each time series
by one single vertex in the graph provides a theoretical framework for
such frequency-domain based analyses.

Secondly, multivariate time series are often obtained by
high-frequency sampling of continuous-time processes such
as EEG-recordings or neural spike trains. Here, our approach
yields a graphical representation of the interrelationships
that does not depend (to some extent) on the sampling frequency
\citep[e.g.,][]{eichlerbrain}. Moreover, many sophisticated models
that have been proposed, for example, for analysing neural activity
do not show a dependence on the past values only at specific lags.
For instance, in the binary time series model discussed in Example
\ref{binarymodel}, the conditional distribution of $X_b(t)$ given
the past history $\ffilt{V}(t-1)$ depends on another process $X_a$
through the past values $X_a(t-1),\ldots,X_a(t-\gamma_b(t))$, where
$\gamma_b(t)$ is the time elapsed since the last event of process $X_b$.
In other words, the number of lagged variables $X_a(t-u)$ on which $X_b(t)$
depends varies over time depending on the past of $X_b$ itself.
Consequently, it seems inappropriate to break down the dependence of
$X_a(t)$ on $\pastx_b(t)$ further into dependencies of $X_a(t)$ on
$X_b(t-u)$ as required by the detailed modelling approach.

\subsection*{\em $m$-separation versus $p$-separation}

The contemporaneous dependence structure of a process $X_V$
can also be described by conditional independencies
of the form
\[
\ffilt{A}(t+1)\indep\ffilt{B}(t+1)\given\ffilt{V}(t),
\]
in which case $X_A$ and $X_B$ are said to be contemporaneously
independent with respect to $\ffilt{V}$. This alternative approach,
which is related to the concept of instantaneous causality by
\citet{granger69}, has been studied by \citep{eichlerpathdiagr} in
the context of weakly stationary processes and linear dependencies.

The most important difference between these two approaches
for defining graphical time series models is that the
corresponding composition and decomposition property
\begin{equation}
\label{condindcondition}
\begin{split}
\ffilt{A}(t+1)\indep &\ffilt{B}(t+1)\given\ffilt{V}(t)\\
&\iff\ffilt{a}(t+1)\indep\ffilt{b}(t+1)\given\ffilt{V}(t)
\quad\forall\,a\in A,\,\forall\,b\in B
\end{split}
\end{equation}
does not follow from condition (S) but requires additional
assumptions similar to condition \eqref{cond block}. Furthermore,
we note that only the first two conditions in Proposition \ref{suffcond}
are sufficient for the above property \eqref{condindcondition}.
Consequently, the class of graphical time series models for which the
pairwise and the block-recursive Granger-causal Markov properties are
equivalent would be smaller under the alternative approach
based on contemporaneous independence. Alternatively, if modelling is to be
based on $m$-separation, one might consider use of an adapted variant of the
connected set Markov property as in \citet{drtonrichardson08} instead
of the pairwise Markov property.

\subsection*{\em Self-loops}

In this paper, we have focused on modelling and analysing the 
interrelationships in multivariate time series. Therefore, we have
not considered the possibility of directed self-loops $v\DE v$,
which could be used to impose additional constraints of the form
$X_B(t+1)\indep\ffilt{B}(t)\given\ffilt{V\without B}(t)$ on a model.
We note that, for a discussion of the dynamic interrelationships among
variables, these self-loops are irrelevant. In fact, it can be shown
that two disjoint sets $A$ and $B$ are $p$-separated given $S$
in a graph with self-loops if and only if they are also $p$-separated
given $S$ in the same graph with all self-loops removed. Similar
statements can be formulated for pointing and extended bi-pointing
paths.

\subsection*{\em Non-stationary time series}

One of our main assumptions has been that the considered multivariate
time series are stationary. This assumption, however, has been made
mainly for the sake of simplicity, and the graphical modelling approach
presented can be extended easily also to the case of non-stationary time
series by requiring that the Granger noncausality and contemporaneous
conditional independence constraints encoded by a graph hold at all
time points in an interval $T\subseteq\znum$, say; in that case, we say
that the time series obeys a Granger-causal Markov property with respect
to the graph over the time interval $T$. This allow us to consider
non-stationary time series models in which the pattern of dependencies
remains fixed whereas the strength of the dependencies may change over
time. An interesting extension would be models where also the
graphical structure changes at certain times. For instance,
\citet{talih05} consider covariance selection models for multivariate
time series where changes in the dependence structure occur at random
times; this approach, however, does not model dynamic dependencies
among the variables. Finally, we note that, despite their practical
relevance, non-stationary models have attracted much less---particularly
theoretical---interest than stationary models due to the involved
inferential problems.

\bigskip

Two important issues have not been addressed in this paper.
Firstly, in many applications there is little prior knowledge about
the causal relationships between the variables, and empirical methods
have to be used to find an appropriate graphical model. This step
of model selection is hampered by the large number of possible models
by which an exhaustive search becomes infeasible even for moderate
dimensions. Therefore, model search strategies are required to lessen
the computational burden.

A second issue, which is related to the problem of model selection,
is the identification of causal effects. It is clear from the definition
of Granger causality that we may conclude from Granger causality to the
existence of a causal effect only if all relevant variables are included
in a study, whereas the omission of important variables can lead
to spurious causalities. However, \citet{hsiao82} noted that such spurious
causalities may vanish if the information set is reduced. In other
words, two processes that both satisfy the pairwise causal Markov
property with respect to a graph $G$ may exhibit different Granger
noncausality relations with respect to partial information sets due
to the presence or absence of spurious causalities. Some concepts
as to how this observation could be exploited for causal inference have been
discussed in \citet{eichlerbrain,ME06latent,eichler:clmps}.

\section*{Acknowledgement}
The author would like to thank two anonymous referees for their
comments and suggestions, which greatly improved the paper.

\begin{appendix}
\section{Conditional independence and stochastic processes}
\label{appendix condind}

Throughout the paper we consider a fixed probability space $(\Omega,\falg,\prob)$.
For any sub-$\sigma$-algebra $\mathscr{H}$ of $\falg$, $\overline{\mathscr{H}}$
denotes the completed $\sigma$-algebra generated by $\mathscr{H}$ and the $\prob$-null
sets in $\falg$. Thus the sets in the completed $\sigma$-algebra $\overline{\mathscr{H}}$
are still measurable sets in $\fclass$. Next, let $\falg_1$, $\falg_2$, and $\falg_3$
be sub-$\sigma$-algebras of $\falg$.
The smallest $\sigma$-algebra generated by $\falg_i\cup\falg_j$ is denoted
as $\falg_i\vee\falg_j$. Then $\falg_1$ and $\falg_2$ are said to be
independent conditionally on $\falg_3$ if
$\mean(X|\falg_2\vee\falg_3)=\mean(X|\falg_3)$ a.s.~for all
real-valued, bounded, $\falg_1$-measurable random variables $X$.
Using the notation of \citet{dawid79} we write
$\falg_1\indep\falg_2\given\falg_3\wrt{\prob}$ or
$\falg_1\indep\falg_2\given\falg_3$ if the reference to $\prob$ is clear.

Let $\falg_i$, $i=1,\ldots,4$ be sub-$\sigma$-algebras of $\falg$.
Then the basic properties of the conditional independence relation are:
\begin{Alist}{CI}
\item
$\falg_1\indep\falg_2\given\falg_3
\,\follows\,\falg_2\indep\falg_1\given\falg_3$
(symmetry)
\item
$\falg_1\indep\falg_2\vee\falg_3\given\falg_4
\,\follows\,\falg_1\indep\falg_2\given\falg_4$ (decomposition)
\item
$\falg_1\indep\falg_2\vee\falg_3\given\falg_4
\,\follows\,\falg_1\indep\falg_2\vee\falg_3\given\falg_3\vee\falg_4$
(weak union)
\item
$\falg_1\indep\falg_2\given\falg_4$ and
$\falg_1\indep\falg_3\given\falg_2\vee\falg_4
\,\follows\,\falg_1\indep\falg_2\vee\falg_3\given\falg_4$ (contraction)
\savecounter{alphcount}
\end{Alist}

In some of the proofs in this paper, we make use of an additional
property,
\begin{Alist}{CI}
\restorecounter{alphcount}
\item
$\falg_{1}\indep \falg_{2}\given \falg_{3}\vee \falg_{4}\text{ and }
\falg_{1}\indep \falg_{3}\given \falg_{2}\vee \falg_{4}
\iff \falg_{1}\indep \falg_{2}\vee \falg_{3}\given \falg_{4}$,
\end{Alist}
which has been called {\em intersection property} by \citet{pearl88}.
Unlike the other basic properties of conditional indepence,
this property does not hold in general. A sufficient and necessary
condition for (CI5) is given by
\begin{equation}
\label{measurable sep}
\overline{\falg_{2}\vee\falg_{4}}
\cap\overline{\falg_{3}\vee\falg_{4}}
=\overline{\falg_{4}}.
\end{equation}
In that case, $\falg_{2}$ and $\falg_{3}$ are said to be measurably
separated conditionally on $\falg_{4}$, denoted by
$\falg_{2}\parallel\falg_{3}\given\falg_{4}\wrt{\prob}$ \citep{fmr90}.
We note that the dependence on $\prob$ is only through the null sets of
$\prob$. For details on conditional measurable separability and its
properties, we refer to Chapter 5.2 of \citet{fmr90}.


If the $\sigma$-algebras $\falg_i$ are generated by random vectors $X_i$
for $i=1,\ldots,4$, in which case we write $\falg_i=\sigma\{X_i\}$, a
sufficient condition for conditional measurable
separability of the $X_i$'s and, thus, of the $\falg_i$'s is that the
probability measure $\prob^{X_1,\ldots,X_4}$ is absolutely continuous with
respect to a product measure $\mu$ and has a positive and continuous density.
However, if each of the $\sigma$-algebras $\falg_i$ is generated by infinitely
many random variables, the condition is obviously no longer valid.
In the following we show that for strictly stationary processes $X_V$ it is
sufficient to assume the existence of a positive and continuous density
for the conditional distribution of $X_V(t+1)$ given its past $\pastx_V(t)$
at the cost of the additional regularity condition (M).

\begin{lemma}
\label{finitemeassep}
Suppose that $X_V$ is a stochastic process such that condition (P) holds,
and let $Y_1,Y_2$ be finite disjoint subsets of
$S(t)=\{X_v(s),s\leq t,v\in V\}$. Then
\begin{equation}
\label{finitemeassep1}
Y_1\parallel Y_2\given\sigma\big\{S(t)\without(Y_1\cup Y_2)\big\}\wrt{\prob},
\end{equation}
where $\sigma\{S(t)\without(Y_1\cup Y_2)\}$ denotes the $\sigma$-algebra
generated by $S(t)\without(Y_1\cup Y_2)$.
\end{lemma}
\begin{proof}
A sufficient condition for \eqref{finitemeassep1}
\citep[Corollary 5.2.11]{fmr90} is the existence of a probability measure
$\prob'$ on $(\Omega,\ffilt{V}(t))$ such that $\prob'$ and
$\prob|_{\ffilt{V}(t)}$, the restriction of $\prob$ on $(\Omega,\ffilt{V}(t))$,
are equivalent (i.e.~have the same null sets) and
\begin{equation}
\label{finitemeassep2}
Y_1\indep Y_2\given\sigma\big\{S(t)\without(Y_1\cup Y_2)\big\}\wrt{\prob'}.
\end{equation}
Take $k\in\nnum$ such that $Y_1\cup Y_2$ and $S(t-k)$ are disjoint, and
let $Z_{j}=X_V(t-j)$ for $j=0,\ldots,k-1$ and $\boldsymbol{Z}_{k}=S(t-k)$.
Noting that by condition (P) the conditional densities
$f_{Z_{jv}|\boldsymbol{Z_{k}}}$ exist and can be derived from the product of
the conditional densities $f_{Z_{j}|Z_{j+1},\ldots,Z_{k-1},\boldsymbol{Z}_{k}}$,
we define the probability kernel $Q(\boldsymbol{z}_{k},A)$ from
$\rnum^{V\times\nnum}$ to $\rnum^{V\times k}$ by
\[
Q(\boldsymbol{z}_{k},A_0\times\cdots\times A_{k-1})
=\int_{A_{k-1}}\!\!\cdots\int_{A_0}
\lprod_{j=0}^{k-1}\lprod_{v\in V}
f_{Z_{jv}|\boldsymbol{Z}_{k}}(z_{jv}|\boldsymbol{z}_k)\,
d\nu(z_{0})\cdots d\nu(z_{k-1}).
\]
Then the probability $\prob'$ on $\big(\Omega,\ffilt{V}(t)\big)$ defined by
\[
\begin{split}
\prob'\big(Z_0\in A_0,\ldots,&Z_{k-1}\in A_{k-1},\boldsymbol{Z}_k\in\boldsymbol{A}_k\big)\\
&=\int_{\boldsymbol{Z}_k^{-1}(\boldsymbol{A}_k)}\int_{A_{k-1}}\!\cdots\int_{A_0}
Q\big(\boldsymbol{Z}_k(\omega),(dz_0,\ldots,dz_{k-1})\big)\,d\prob(\omega)
\end{split}
\]
is equivalent to $\prob|_{\ffilt{V}(t)}$. Furthermore, the random variables
$Z_{jv}$ with $j=0,\ldots,k-1$ and $v\in V$ are mutually independent
conditionally on $\boldsymbol{Z}_k$ under $\prob'$, which implies
\eqref{finitemeassep2} and hence \eqref{finitemeassep1}.
\end{proof}

The next result shows that this conditional measurable separability can also
be extended to $\sigma$-algebras $\ffilt{A}(t)$ generated by the pasts
$\pastx_A(t)$ provided the process $X_V$ is conditionally mixing (in the sense
of condition (M)).

\begin{proposition}
\label{meassep-lemma}
Suppose that $X_V$ is a stochastic process such that conditions (M) and (P)
hold. Then $\ffilt{A}(t)$ and $\ffilt{B}(t)$ are measurably separated
conditionally on $\ffilt{{V\without(A\cup B)}}(t)$ for all disjoint
subsets $A$ and $B$ of $V$ and all $t\in\znum$.
\end{proposition}
\begin{proof}
Let $A$ and $B$ be disjoint subsets of $V$. We have to show that
$\ffilt{A}(t)$, $\ffilt{B}(t)$, and $\ffilt{{V\without(A\cup B)}}(t)$
satisfy \eqref{measurable sep} and hence that
\begin{equation}
\label{measseplemma1}
\overline{\ffilt{{V\without B}}(t)}\cap\overline{\ffilt{{V\without A}}(t)}
=\overline{\ffilt{{V\without(A\cup B)}}(t)}
\end{equation}
for all $t\in\znum$. From Lemma \ref{finitemeassep}, it follows that,
for all $t\in\znum$ and $k\in \nnum$, the $\sigma$-algebras
$\sigma\{X_A(t),\ldots,X_A(t-k+1)\}$ and $\sigma\{X_B(t),\ldots,X_B(t-k+1)\}$
are measurably separable conditional on
$\ffilt{{V\without(A\cup B)}}(t)\vee\ffilt{V}(t-k)$.
Accordingly, we have by the definition of conditionally measurable
separability
\[
\overline{\ffilt{{V\without B}}(t)\vee\ffilt{V}(t-k)}\cap
\overline{\ffilt{{V\without A}}(t)\vee\ffilt{V}(t-k)}
=\overline{\ffilt{{V\without(A\cup B)}}(t)\vee\ffilt{V}(t-k)}
\]
for all $t\in\znum$ and $k\in\nnum$. Since the $\sigma$-algebras on both
sides are monotonically decreasing as $k$ increases, this yields 
for $k\to\infty$
\begin{align*}
\mathop{\textstyle\bigcap}_{k>0}\big[
\overline{\ffilt{{V\without B}}(t)\vee\ffilt{V}(t-k)}\cap
&\overline{\ffilt{{V\without A}}(t)\vee\ffilt{V}(t-k)}\big]
=\mathop{\textstyle\bigcap}_{k>0}
\overline{\ffilt{{V\without(A\cup B)}}(t)\vee\ffilt{V}(t-k)}
\end{align*}
for all $t\in\znum$. Since by condition (M)
\[
\mathop{\textstyle\bigcap}_{k>0}\big[\overline{\ffilt{{S}}(t)}
\vee\overline{\ffilt{V}(t-k)}\big]
=\overline{\ffilt{{S}}(t)}
\]
for all subsets $S$ of $V$, this establishes \eqref{measseplemma1}.
\end{proof}

\begin{proof}[Proof of Proposition \ref{intersectprop}]
The result follows directly from Lemma \ref{finitemeassep} and
Proposition \ref{meassep-lemma}.
\end{proof}

\section{Graphical terminology}

We firstly recall some basic graphical definitions used in this paper.
In a graph $G=(V,E)$, if there is a directed edge $a\DE b$, we say that
$a$ is a parent of $b$ and $b$ is a child of $a$; similarly, if there is
an undirected line $a\UE b$, the vertices $a$ and $b$ are called neighbours.
The sets of parents, children and neighbours of a vertex $a$ are
denoted as $\parent{a}$, $\child{a}$, and $\neighbour{a}$, respectively.
Furthermore, for $A\subseteq V$, let $\parent{A}=\cup_{a\in A}\parent{a}\without A$ be the set of all parents
of vertices in $A$ that are not themselves in $A$, and let $\child{A}$ and
$\neighbour{A}$ be defined similarly.

Next, as in \citet{frydenberg}, a vertex $b$ is said to be an {\em ancestor}
of $a$ if either $b=a$ or there exists a directed path $b\DE\cdots\DE a$ in
$G$. The set of all ancestors of elements in $A$ is denoted by $\ancestor{A}$.
Notice that this definition differs from the one given in \citet{SL96}.
A subset $A$ is called an {\em ancestral set} if it contains
all its ancestors, that is, $\ancestor{A}=A$.

Finally, let $G=(V,E)$ and $G'=(V',E')$ be mixed graphs. Then $G'$ is a
{\em subgraph} of $G$ if $V'\subseteq V$ and $E'\subseteq E$.
If $A$ is a subset of $V$ it induces the
subgraph $G_A=(A,E_A)$ where $E_A$ contains all edges $e\in E$
that have both endpoints in $A$.

In the remainder of this section, we prove some auxiliarly results
that allow us to relate separation statements in the full graph to
separation statement in so-called marginal graphs, which
basically reflect the dynamic dependencies in appropriate subprocesses
(see Lemma \ref{marginal bc}).

\begin{definition}[Marginal graph]
Let $G=(V,E)$ be a mixed graph and let $A$ be an ancestral subset of $V$.
Then the {\em marginal graph} $G_{[A]}=(A,E_{[A]})$ induced by $A$
is obtained from the induced subgraph $G_{A}$ by insertion of additional
undirected edges $a\UE b$ whenever there exists an undirected path between
$a$ and $b$ in $G$ that does not intersect $\ancestor{A}\without\{a,b\}$.
\end{definition}

\begin{lemma}
\label{psepancestral}
Let $G=(V,E)$ be a mixed graph and $A$, $B$, $S$ disjoint subsets of $V$. 
Then $A$ and $B$ are $p$-separated given $S$ in $G$ if and only if
$A$ and $B$ are $p$-separated given $S$ in $\GAG{A\cup B\cup S}$.
\end{lemma}
\begin{proof}
To show necessity, let $\pi=\path{e_1,\ldots,e_n}$ be a $p$-connecting
path between $A$ and $B$ given $S$ in $\GAG{A\cup B\cup S}$.
If all edges of $\pi$ are edges in $G$, $\pi$ is also $p$-connecting given $S$ in $G$. Thus,
we may assume that there exist edges in $\pi$, $e_{j_1},\ldots,e_{j_m}$ say,
that do not occur in $G$. These edges $e_{j_k}$ are necessarily undirected
since all directed edges in $\GAG{A\cup B\cup S}$ also occur in $G$.
Let $e_{j_k}=v_{j_k}\UE v_{j_k+1}$. Then by definition of the marginal
graph there exists an undirected path $\phi_{j_k}$ between
$v_{j_k-1}$ and $v_{j_k}$ which bypasses
$\ancestor{A\cup B\cup S}\without\{v_{j_k-1},v_{j_k}\}$
and therefore is $p$-connecting given $S$.
Replacing all edges $e_{j_k}$ in $\pi$ by the corresponding paths
$\phi_{j_k}$ we obtain a new path $\pi'$ which connects $A$ and $B$ in $G$.
This path $\pi'$ is also $p$-connecting given $S$ since the replacement
of $e_{j_k}$ by the undirected and $p$-connecting path $\phi_{j_k}$
does not change the $p$-collider resp.~$p$-noncollider status of the
points $v_{j_k-1}$ and $v_{j_k}$.

Conversely for sufficiency, let $\pi=\path{e_1,\ldots,e_n}$ be a
$p$-connecting path between $A$ and $B$ given $S$ in $G$.
Then all edges in $\pi$ with both endpoints in $\ancestor{A\cup B\cup S}$
also occur in $\GAG{A\cup B\cup S}$ since $G_{\ancestor{A\cup B\cup S}}$ is a subgraph of
$\GAG{A\cup B\cup S}$.
We firstly show that the endpoints of any directed edge $e_j$
in $\pi$ are in $\ancestor{A\cup B\cup S}$. Let $e_j=v_j\DE v_{j+1}$
(the case $e_j=v_j\LDE v_{j+1}$ is treated similarly).
Then there exists a directed subpath $\path{e_j,\ldots,e_{j+r}}$
of maximal length such that either $v_{j+r}$ is an endpoint of
$\pi$ and, thus, in $A\cup B$ or
$e_{j+r+1}$ is of the form $v_{j+r}\UE v_{j+r+1}$ or $v_{j+r}\LDE v_{j+r+1}$.
In the latter case $v_{j+r}$ is a $p$-collider and, thus, in $S$
since $\pi$ is $p$-connecting given $S$.
It follows that $v_{j}$ and $v_{j+1}$ are both in $\ancestor{A\cup B\cup S}$.

Next, if $e_{j}$ is an edge in $\pi$ that does not occur in
$\GAG{A\cup B\cup S}$, at least one of its endpoints $v_{j-1}$ and
$v_j$ is not in $\ancestor{A\cup B\cup S}$.
Thus, there exists an undirected subpath
$\psi_{i,k}=\path{e_i,\ldots,e_{k}}$ with $i\leq j\leq k$
such that $v_{i-1},v_{k}\in\ancestor{A\cup B\cup S}$
but all intermediate points are not in $\ancestor{A\cup B\cup S}$.
In other words, $v_{i-1}$ and $v_{k}$ are not separated by
$\ancestor{A\cup B\cup S}\without\{v_{j-1},v_{k}\}$ in $G$
which implies the presence of the undirected edge
$f_{i,k}=v_{i-1}-v_{k}$ in $\GAG{A\cup B\cup S}$. Replacing all
undirected subpaths $\phi_{i,k}$ with intermediate points not in
$\ancestor{A\cup B\cup S}$ by the corresponding edge $f_{i,k}$,
we obtain a path between $A$ and $B$ in $\GAG{A\cup B\cup S}$
which still has all its $p$-collider in $S$ and all its $p$-noncolliders
outside $S$ and therefore is $p$-connecting given $S$.
\end{proof}

The following lemma is an adapted version of Proposition 2 in
\citet{koster99}. The proof is considerably shorter due to the
fact that we allow paths to be self-intersecting.

\begin{lemma}
\label{fullpsep}
Let $A$, $B$, $S$ be disjoint subsets of $V$. Then $A$ and $B$ are
$p$-separated given $S$ in $G_{[\ancestor{A\cup B\cup S}]}$
if and only there exist subsets $A'$ and $B'$
such that $A\subseteq A'$, $B\subseteq B'$,
$A'\cup B'\cup S=\ancestor{A\cup B\cup S}$
and
\[
A'\sepp B'\given S\wrt{G_{[\ancestor{A\cup B\cup S}]}}.
\]
\end{lemma}
\begin{proof}
By Lemma \ref{psepancestral} we may assume that $V=\ancestor{A\cup B\cup S}$.
Let $A'$ be the subset of vertices $v\in V\without(B\cup S)$ such that
$v\sepp B\given S\wrt{G}$, and set $B'=V\without(A'\cup S)$.
Then $A'$ and $B$ are obviously $p$-separated given $S$. Thus, we have
to show that $a$ and $b'$ are $p$-separated given $S$ whenever $a\in A'$ and
$b'\in B'\without B$. Suppose to the contrary that there exists a
$p$-connecting path $\pi$ between some $a\in A'$ and $b'\in B'\without B$.
Since $A'$ contains all vertices in $V\without(B\cup S)$ that are
$p$-separated from $B$ given $S$, there exists a $p$-connecting path
$\pi'$ between $b'$ and some $b\in B$. Furthermore, since
$b'\in\ancestor{A\cup B\cup S}\without(A\cup B\cup S)$ there exists
some vertex $u\in A\cup B\cup S$ and a directed path
$\omega=b'\DE\cdots\DE u$ with no intermediate points in $A\cup B\cup S$.
Denoting by $\rev{\omega}$ the reverse path of $\omega$, that is,
$\rev{\omega}=u\LDE\cdots\LDE b'$, we may compose a path $\phi$ between
$A$ and $B$ by
\begin{romanlist}
\item
$\phi=\path{\rev{\omega},\pi'}$ if $u\in A$,
\item
$\phi=\path{\pi,\omega}$ if $u\in B$, and
\item
$\phi=\path{\pi,\omega,\rev{\omega},\pi'}$ if $u\in S$.
\end{romanlist}
We note that the directed path $\omega$ is $p$-connecting given $S$
since it has no intermediate points in $S$. Furthermore, $b'\notin S$
is a $p$-noncollider on $\phi$ in each of these cases and $v\in S$ is a
$p$-collider on $\phi$ in case (iii). Hence $\phi$ is a $p$-connecting
path between $A$ and $B$ given $S$ which contradicts our assumption.

The opposite implication is obvious because of the elementwise
definition of $p$-separation.
\end{proof}

Because of Lemmas \ref{psepancestral} and \ref{fullpsep}, it is often
sufficient in the proofs to consider only the case of $A\sepp B\given S$
with $S=V\without(A\cup B)$. In this case, $p$-separation can be
characterized in terms of pure-collider paths---paths on which every
intermediate node is a collider---or in terms of local configurations.

\begin{lemma}
\label{augmentationrules}
Let $G$ be a mixed graph and let $A$ and $B$ be two disjoint subsets
of $V$. Then the following statements are equivalent:
\begin{romanlist}
\item
$A\sepp B\given V\without(A\cup B)$;
\item
$A$ and $B$ are not connected by a pure-collider path;
\item
$(A\cup\child{A})\cap(B\cup\child{B})=\varnothing$ and
$\neighbour{A\cup\child{A}}\cap(B\cup\child{B})=\varnothing$.
\end{romanlist}
\end{lemma}

Note that the second part of condition (iii) states that no two vertices
$a\in A\cup\child{A}$ and $b\in B\cup\child{B}$ are adjacent; the condition
thus is also symmetric in $A$ and $B$.

\begin{proof}
This observation follows directly from the definition of
$p$-separation and pure-collider paths.
\end{proof}

\section{Proofs}

\begin{proof}[Proof of Theorem \ref{basicMP}]
Setting $A=\{a\}$ in (BC), we obtain (LC). Conversely,
since $\parent{a}\cup\{a\}\subseteq\parent{A}\cup A$, we have
by (LC) together with (CI2) and (CI3)
\[
X_{V\without\parent{A}\cup A}\noncausal X_{a}\quad\forall a\in A,
\]
which, under condition \eqref{cond block}, implies the first part of (BC).
The second part is proved similarly.

To see that (LC) and (PC) are equivalent, we note that, under condition
(S), the intersection property leads to the following
composition and decomposition property for Granger noncausality
relations:
\begin{equation}
\label{comp/decomp}
X_A\noncausal X_B\wrt{\ffilt{V}}\iff
X_{a}\noncausal X_B\wrt{\ffilt{V}}\quad\forall a\in A.
\end{equation}
Similarly, we have for contemporaneous conditional independence relations
\begin{equation}
\label{contemp block}
X_A\noncorr X_B\wrt{\ffilt{V}}\iff
X_a\noncorr X_b\wrt{\ffilt{V}}\quad\forall\,a\in A,\,\forall\,b\in B.
\end{equation}
Taking $A=V\without(B\cup\parent{B})$ in \eqref{comp/decomp} and
$A=V\without(B\cup\parent{B})$ in \eqref{contemp block}, we find that
the pairwise and the local Granger-causal Markov properties are equivalent.
\end{proof}

\begin{proof}[Proof of Proposition \ref{suffcond}]
By Theorem \ref{basicMP}, it suffices to show that each of the three
conditions (i), (ii), and (iii) implies
\begin{equation}
\label{suffcondeq}
X_A\noncausal X_b\wrt{\ffilt{V}}\quad\forall\,b\in B\follows
X_A\noncausal X_B\wrt{\ffilt{V}}
\end{equation}
for any two disjoint subsets $A,B\subseteq V$.

For the first case, let $H$ be the Hilbert space of all square
integrable random variables on $(\Omega,\falg,\prob)$. Furthermore,
for $U\subseteq V$, let $H_U(t)$ be the closed subspace spanned
by $\{X_u(s),u\in U,s\leq t\}$ and let $H^\perp_U(t)$ be its
orthogonal complement. Then we have for any $Y\in H^\perp_{V\without A}(t)$
\[
\cov\big(X_B(t+1),Y\big)=0
\iff\cov\big(X_b(t+1),Y\big)=0\quad\forall\,b\in B,
\]
which for a Gaussian process implies \eqref{suffcondeq}.

Next, suppose that condition (ii) holds and that $X_A$ is Granger-noncausal
for $X_b$ with respect to $\ffilt{V}$ for all $b\in B$. Then, the conditional
distribution $\prob^{X_B(t+1)|\pastx_{V}(t)}$ satisfies
\[
\prob^{X_B(t+1)|\pastx_{V}(t)}
=\otimes_{b\in B}\prob^{X_b(t+1)|\pastx_{V}(t)}
=\otimes_{b\in B}\prob^{X_b(t+1)|\pastx_{V\without A}(t)}
\]
and, thus, is $\ffilt{V\without A}(t)$-measurable,
which proves \eqref{suffcondeq}.

Finally, if condition (iii) holds, we have
\[
X_B(t+1)-\mean\big[X_B(t+1)\given\ffilt{V}(t)\big]\indep
\ffilt{A}(t)\given\ffilt{V\without A}(t).
\]
Since the left hand side of \eqref{suffcondeq} implies that
$\mean\big[X_B(t+1)\given\ffilt{V}(t)\big]$ is
$\ffilt{V\without A}(t)$-measurable,
we obtain $X_B(t+1)\indep\ffilt{A}(t)\given\ffilt{V\without A}(t)$,
which completes the proof.
\end{proof}

For the proof of the equivalence of the block-recursive and the global
Granger-causal Markov property, it will be convenient to restrict ourselves to
mixed graphs for ancestral subsets.
Due to the additional undirected edges inserted into the marginal graph
$\GAG{A}$, the subprocess $X_{\ancestor{A}}$ satisfies the pairwise
Granger-causal Markov property with respect to $\GAG{A}$ if $X_V$ did so
with respect to $G$. The following lemma shows that the same inheritance
property also holds for the block-recursive Granger-causal Markov property.

\begin{lemma}
\label{marginal bc}
Suppose that $X_V$ satisfies the block-recursive Granger-causal Markov property
with respect to the mixed graph $G$, and let $U\subseteq V$. Then the
subprocess $X_{\ancestor{U}}$ satisfies the block-recursive Granger-causal
Markov property with respect to the marginal ancestral graph $\GAG{U}$.
\end{lemma}

\begin{proof}
Let $H=\GAG{U}$ and let $A$ be a subset of $\ancestor{U}$. We first note that,
since $\ancestor{U}$ is an ancestral set and, thus, contains the parents
of all its subsets $A$, the parents of $A$ in both graphs are the same,
that is, $P=\parentin{G}{A}=\parentin{H}{A}$. By the block-recursive
Granger-causal Markov property of $X_V$ with respect to $G$,
$X_{V\without(P\cup A)}$ does not Granger-cause $X_A$
with respect to $\ffilt{V}$, which by (CI2) implies that
$X_{\ancestor{U}\without(P\cup A)}$ is Granger-noncausal for $X_A$ with
respect to the smaller filtration $\ffilt{\ancestor{U}}$ as required by the
block-recursive Granger-causal Markov property of
$X_{\ancestor{U}}$ with respect to $H$.

Next, let $N=\neighbourin{H}{A}$. Then $A$ and $\ancestor{U}\without (N\cup A)$
are separated by $N$ in $H\undirect$, that is, $a$ and $b$ are not
adjacent in the undirected subgraph $H\undirect$ whenever $a\in A$
and $b\in\ancestor{U}\without(N\cup A)$.
By definition of $H$, this implies that $A$ and $\ancestor{U}\without(N\cup A)$
are separated by $N$ in $G\undirect$. By the block-recursive
Granger-causal Markov property, it follows that
\[
\ffilt{A}(t+1)\indep\ffilt{\ancestor{U}\without(N\cup A)}(t+1)\given
\ffilt{V}(t)\vee\ffilt{N}(t+1)
\]
and, with (CI2) and (CI3),
\[
\ffilt{A}(t+1)\indep\ffilt{V\without\ancestor{U}}(t)\given
\ffilt{\ancestor{U}\cup N}(t).
\] 
Combining these two relations by using (CI2) to (CI4), we find that
$X_{\ancestor{U}\without(N\cup A)}$ and $X_A$ are contemporaneously
conditionally independent with respect to $\ffilt{\ancestor{U}}$ as
required by the block-recursive Granger-causal Markov property of
$X_{\ancestor{U}}$ with respect to the graph $H$.
\end{proof}

\begin{figure}
\centerline{\includegraphics[width=9cm]{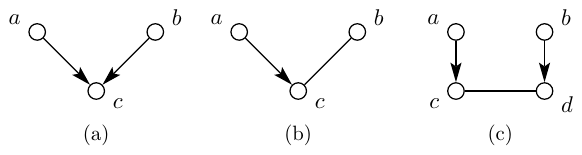}}
\caption{Pure-collider paths between two vertices $a$ and $b$.}
\label{moralrules}
\end{figure}

\begin{proof}[Proof of Lemma \ref{globalMPpast}]
For notational convenience, we may assume in view of Lemma \ref{marginal bc}
that $\ancestor{A\cup B\cup S}=V$ and, thus, $\GAG{A\cup B\cup S}=G$.
Furthermore, Lemma \ref{psepancestral} implies that, if $A\sepp B\given S$
in the graph $G$, there exists a partition $(A^*,B^*,S)$ of $V$ such that
$A\subseteq A^*$, $B\subseteq B^*$, and $A^*\sepp B^*\given S$.
Thus, without loss of generality, we may assume that
$S=V\without(A\cup B)$.

With these simplifications, it suffices to show that
$A\sepp B\given V\without(A\cup B)$ implies
\begin{equation}
\label{pr:MPpast1}
\ffilt{X_A}(t)\indep\ffilt{X_B}(t)\given\ffilt{{V\without(A\cup B)}(t)}(t)
\end{equation}
for all $t\in\znum$. To this end, we firstly show that
\begin{equation}
\label{globalMPpast0}
\ffilt{A}(t)\indep\ffilt{B}(t)\given
\ffilt{V\without(A\cup B)}(t)\vee\ffilt{A\cup B}(t-k)
\end{equation}
for all $t\in\znum$ and $k\in\nnum$.

We proceed by induction on $k$. For $k=1$, we obtain \eqref{globalMPpast0}
immediately from the block-recursive Granger-causal Markov
property noting that $B\subseteq V\without(A\cup\neighbour{A})$.
For the induction step $k\to k+1$ assume that
\begin{equation}
\label{inductassump}
\ffilt{A}(t)\indep\ffilt{B}(t)\given
\ffilt{V\without(A\cup B)}(t)\vee\ffilt{A\cup B}(t-k)
\end{equation}
for all $t\in\znum$. Let $C_A=A\cup\child{A}$. Then, since by the
block-recursive Granger-causal Markov property $X_A$ is Granger-noncausal
for $X_{V\without C_A}$ with respect to $\ffilt{V}$, we have
\[
\ffilt{A}(t)\indep\ffilt{V\without C_A}(t+1)
\given\ffilt{V\without A}(t)\vee\ffilt{A\cup B}(t-k)
\]
and further with \eqref{inductassump} and (CI4)
\[
\ffilt{A}(t)\indep\ffilt{B}(t)\vee\ffilt{V\without C_A}(t+1)
\given\ffilt{V\without(A\cup B)}(t)\vee\ffilt{V}(t-k).
\]
With $N_A=\neighbour{A\cup\child{A}}=\neighbour{C_A}$,
we obtain by (CI2) and (CI3)
\begin{equation}
\label{pr:MPpast2}
\ffilt{A}(t)\indep\ffilt{B}(t)\vee\ffilt{V\without(C_A\cup N_A)}(t+1)
\given\ffilt{N_A}(t+1)\vee\ffilt{V\without(A\cup B)}(t)\vee\ffilt{V}(t-k).
\end{equation}

Next, we note that by Lemma \ref{augmentationrules}
$B\cup\child{B}\subseteq V\without(C_A\cup N_A)$
and thus
\[
\ffilt{C_A}(t+1)\indep\ffilt{B}(t)\given
\ffilt{N_A}(t+1)\vee\ffilt{V\without B}(t).
\]
Furthermore, $X_{C_A}$ and
$X_{V\without(C_A\cup N_A)}$ are contemporaneously conditionally
independent and thus
\[
\ffilt{C_A}(t+1)\indep\ffilt{V\without (C_A\cup N_A)}(t+1)
\given\ffilt{N_A}(t+1)\vee\ffilt{V}(t).
\]
Together with the previous relation, we obtain by (CI4)
\[
\ffilt{C_A}(t+1)\indep\ffilt{B}(t)\vee\ffilt{V\without(C_A\cup N_A)}(t+1)
\given\ffilt{N_A}(t+1)\vee\ffilt{V\without B}(t).
\]
By \eqref{pr:MPpast2} together with properties (CI2), (CI3), and (CI5),
this yields
\[
\ffilt{A}(t)\vee\ffilt{C_A}(t+1)
\indep\ffilt{B}(t)\vee\ffilt{V\without(C_A\cup N_A)}(t+1)
\given\ffilt{N_A}(t+1)\vee\ffilt{V\without(A\cup B)}(t)\vee\ffilt{V}(t-k).
\]
Since this relation holds for all $t\in\znum$, we have
by (CI2) and (CI3)
\[
\ffilt{A}(t)\indep\ffilt{B}(t)\given
\ffilt{V\without(A\cup B)}(t)\vee\ffilt{A\cup B}(t-k-1),
\]
which completes the induction step.

To show that \eqref{globalMPpast0} entails \eqref{pr:MPpast1},
we note that for $k\to\infty$ \eqref{globalMPpast0} yields
\[
\ffilt{A}(t)\indep\ffilt{B}(t)
\given
\mathop{\textstyle\bigcap}_{k>0}\big[
\ffilt{{V\without(A\cup B)}}(t)\vee\ffilt{{A\cup B}}(t-k)\big]
\]
for all $t\in\znum$. As in the proof of Proposition \ref{meassep-lemma},
it follows that
\[
\mathop{\textstyle\bigcap}_{k>0}\big[
\overline{\ffilt{{V\without(A\cup B)}}(t)}\vee
\overline{\ffilt{{A\cup B}}(t-k)}\big]
=\overline{\ffilt{{V\without(A\cup B)}}(t)},
\]
which concludes the proof of \eqref{pr:MPpast1}.
\end{proof}

\begin{proof}[Proof of Theorem \ref{global AMP}]
Suppose that $A$, $B$, and $S$ are disjoint subsets of $V$ such that
$A\sepp B\given S$. Let $\xi$ be any $\ffilt{A}(\infty)$
measurable random variable with $\mean|\xi|<\infty$, where
$\ffilt{A}(\infty)=\mathop{\vee}_{t\in\znum}\ffilt{A}(t)$
denotes the $\sigma$-algebra generated by $X_A$.
Then $\xi(t)=\mean\big(\xi|\ffilt{A}(t)\big)$ is a martingale and
converges to $\xi$ in $L^1$ as $t$ tends to infinity.
Thus, we obtain on the one hand, as $t\to\infty$,
\begin{equation}
\label{condmeanconv1}
\mean\big(\xi(t)|\ffilt{S\cup B}(t)\big)\to
\mean\big(\xi|\ffilt{S\cup B}(\infty)\big)\quad\text{in $L^1$.}
\end{equation}
On the other hand, since $\xi(t)\indep\ffilt{B}(t)\given\ffilt{S}(t)$
by Lemma \ref{globalMPpast}, we have, as $t\to\infty$,
\begin{equation}
\label{condmeanconv2}
\mean\big(\xi(t)|\ffilt{S\cup B}(t)\big)
=\mean\big(\xi(t)|\ffilt{S}(t)\big)
\to\mean\big(\xi|\ffilt{S}(\infty)\big)\quad\text{in $L^1$.}
\end{equation}
Since the limits in \eqref{condmeanconv1} and \eqref{condmeanconv2} must be
equal in $L^1$ and, thus, also almost surely, this proves that
$\ffilt{A}(\infty)\indep\ffilt{B}(\infty)\given\ffilt{S}(\infty)$.
\end{proof}

\begin{proof}[Proof of Theorem \ref{global-markov}]
For the proof of the first part of the global Granger-causal Markov property,
let $A$ and $B$ be subsets such that all $B$-pointing paths between $A$ and $B$
are $p$-blocked given $B\cup S$. We note that each $B$-pointing path $\pi$
is of the form $\pi=\path{\tilde\pi,e}$, where $e$ is a directed edge
$u\DE b$ for some $b\in B$. Thus, $\pi$ is $p$-blocked given $B\cup S$
if and only if $u\in B\cup S$ or $\tilde\pi$ is $p$-blocked given
$B\cup S$. Therefore, if all $B$-pointing paths between $A$ and $B$
are $p$-blocked given $B\cup S$, then $A$ and $\parent{B}\without(B\cup S)$
are $p$-separated given $B\cup S$ and we obtain by Lemma \ref{globalMPpast}
\[
\ffilt{\parent{B}\without(B\cup S)}(t)\indep
\ffilt{A}(t)\given\ffilt{B\cup S}(t).
\]
Since, in particular, every edge $a\DE b$ for some $a\in A$ and $b\in B$
is $p$-connecting, it follows that $A$ and $\parent{B}$ are disjoint.
Thus, we get by the block-recursive Granger-causal Markov property
\[
\ffilt{B}(t+1)\indep\ffilt{A}(t)
\given\ffilt{\parent{B}\cup S\cup B}(t).
\]
Applying the contraction property to this and the previous relation,
we find that $X_A$ is Granger-noncausal for $X_B$ with respect to
$\ffilt{A\cup B\cup S}$.

For the proof of the second part, let $U=A\cup B\cup S$ and assume that
every extended bi-pointing path between $A$ and $B$ is $p$-blocked
given $U$. This includes in particular all bi-pointing paths $\pi$ between
$a\in A$ and $b\in B$, which are of the form $\pi=\path{e_1,\tilde\pi,e_n}$,
where $e_1$ and $e_n$ are directed edges $a\LDE p_a$ and $p_b\DE b$,
respectively (Fig.~\ref{extbidirect}\,a). Then $\pi$ is $p$-blocked given
$U$ if and only if $p_a\in U$, $p_b\in U$, or $\tilde\pi$ is $p$-blocked
given $U$. This implies that, if all bi-pointing paths between $A$ and $B$ are
$p$-blocked given $U$, $\parent{A}\without U$ and $\parent{B}\without U$
are $p$-separated given $U$.

\begin{figure}
\includegraphics[width=11cm]{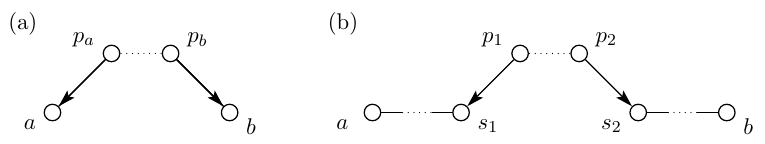}
\caption{(a) bi-pointing path; (b) extended bi-pointing path.}
\label{extbidirect}
\end{figure}

Next, we seek to find subsets $S_A$ and $S_B$ of $S$ such that all extended
bi-pointing paths between the enlarged sets $A\cup S_A$ and $B\cup S_B$ are
still $p$-blocked given $U$. Then, by the same argument as
above, $\parent{A\cup S_A}\without U$ and $\parent{B\cup S_B}\without U$
are $p$-separated given $U$. As an example, consider the extended bi-pointing
path in Fig.~\ref{extbidirect}(b) and suppose that $s_1$ and $s_2$ are
linked to $a$ and $b$, respectively, by undirected paths that are
$p$-connecting given $S$. Then the depicted extended bi-pointing path
is $p$-blocked given $U$ if and only if $p_1$ and $p_2$ are $p$-separated
given $U$.

For a formal definition of the sets $S_A$ and $S_B$, we first set
$S_0=\{s\in S|\parent{s}\subseteq U\}$, which in particular includes all
$s\in S$ that have no parents. Then adding any vertex in $S_0$ to either
$S_A$ or $S_B$ to either $A\cup S_A$ or $B\cup S_B$ will not increase the
sets $\parent{A\cup S_A}\without U$ or $\parent{B\cup S_B}\without U$.Therefore, we set
For a formal argument, we need to define the sets $S_A$ and $S_B$ slightly
differently. More precisely, let $S_0=\{s\in S|\parent{s}\subseteq U\}$,
which in particular includes all $s\in S$ that have no parents. Furthermore,
let $S_A$ be the set of all $s\in S\without S_0$ such that every
extended bi-pointing path between $s$ and $B$ is $p$-blocked given $U$
and set $S_B=S\without(S_0\cup S_A)$. Notice that for all $s\in S_B$
there exists an extended bi-pointing path between $s$ and $B$ that
is $p$-connecting given $U$. We show that every extended bi-pointing path
between $A\cup S_A$ and $B\cup S_B$ is $p$-blocked given $U$. Since all
extended bi-pointing paths between $A\cup S_A$ and $B$ must be $p$-blocked by
assumption on $A$ and $B$ or by definition of $S_A$, we only have to show
that all extended bi-pointing paths between $A\cup S_A$ and $S_B$ are
$p$-blocked given $U$. Suppose to the contrary that $\pi$ is an extended
bi-pointing path between $A\cup S_A$ and $s\in S_B$ that is $p$-connecting
given $U$. Then, as mentioned above, there exists a $p$-connecting extended
bi-pointing path $\pi_s$ between $s$ and $B$. If $s$ is a $p$-collider on
the composed extended bi-pointing path $\tilde\pi=\path{\pi,\pi_s}$ then
$\tilde\pi$ is $p$-connecting given $U$ contradicting the assumption about
$A$ and $B$. Otherwise, if $s$ is a $p$-noncollider, the two adjacent edges
must be undirected (i.e.~$\UE s\UE$) because extended bi-pointing paths
never have a tail at either endpoint. Since $s\notin S_0$ there exists a
path $\check\pi=\path{\pi,s\LDE v\DE s,\pi_s}$ with $v\in\parent{s}\without U$.
The two instances of $s$ on $\check{\pi}$ that are adjacent to $v$ are
$p$-colliders and $\check\pi$ thus is $p$-connecting given $U$. Since
$\check\pi$ is composed of extended bi-pointing paths, it is itself an
extended bi-pointing path. Thus, by definition of
$S_A$, $\check\pi$ must have endpoints in $A$ and $B$, which contradicts again
the assumption about $A$ and $B$.

Since in particular all bi-pointing paths between $A\cup S_A$ and $B\cup S_B$
are $p$-blocked given $U$, we have
\begin{equation}
\label{equiv-0}
\parent{A\cup S_A}\without U\sepp\parent{B\cup S_B}\without U\given U.
\end{equation}
Thus, we obtain by Lemma \ref{globalMPpast}
\begin{equation}
\label{equiv-1}
\ffilt{\parent{A\cup S_A}\without U}(t)\indep
\ffilt{\parent{B\cup S_B}\without U}(t)\given \ffilt{U}(t).
\end{equation}
It also follows from \eqref{equiv-0} that the sets $\parent{A\cup S_A}\without U$
and $\parent{B\cup S_B}$ are disjoint and thus
$\parent{A\cup S_A}\without U\subseteq V\without\parent{B\cup S_B}$,
Noting furthermore that $\parent{S_0}\subseteq U$ by definition of $S_0$,
we obtain from the block-recursive Granger-causal Markov property
\begin{equation}
\label{equiv-2}
\ffilt{B\cup S_B\cup S_0}(t+1)\indep\ffilt{\parent{A\cup S_A}\without U}(t)
\given\ffilt{U\cup \parent{B\cup S_B}}(t).
\end{equation}
Together with \eqref{equiv-1} this yields
\begin{equation}
\label{equiv-3}
\ffilt{B\cup S_B\cup S_0}(t+1)\indep\ffilt{\parent{A\cup S_A}\without U}(t)
\given\ffilt{U}(t).
\end{equation}

Moreover, since undirected paths are special cases of extended bi-pointing
paths, we find that every undirected path between $A\cup S_A$ and $B\cup S_B$
intersects $S_0$. Then, by a standard argument of graph theory
\citep[e.g.,][Lemma 3.3.3]{whittaker90}, there exists a partition $(A^*,B^*,S_0)$ of $V$
such that $A\cup S_A\subseteq A^*$, $B\cup S_B\subseteq B^*$, and
every undirected path between $A^*$ and $B^*$ intersects $S_0$; in particular,
this implies $\neighbour{A\cup S_A}\subseteq S_0$. Thus, we obtain by the
block-recursive Granger-causal Markov property
\[
\ffilt{A\cup S_A}(t+1)\indep\ffilt{B\cup S_B}(t+1)\given\ffilt{V}(t)\vee
\ffilt{S_0}(t+1).
\]
Together with
\[
\ffilt{A\cup S_A\cup S_0}(t+1)\indep
\ffilt{V\without(U\cup\parent{A\cup S_A})}(t)
\given\ffilt{U\cup\parent{A\cup S_A}}(t),
\]
which, by $\parent{S_0}\cup(A\cup S_A\cup S_0)\subseteq U$, also follows
from the block-recursive Granger-causal Markov property, this implies
\begin{equation}
\label{equiv-4}
\ffilt{A\cup S_A}(t+1)\indep\ffilt{B\cup S_B}(t+1)\given
\ffilt{U\cup \parent{A\cup S_A}}(t),X_{S_0}(t+1).
\end{equation}
Applying (CI4) to \eqref{equiv-3} and \eqref{equiv-4},
we finally obtain
\[
\ffilt{A\cup S_A}(t+1)\indep\ffilt{B\cup S_B}(t+1)\given\ffilt{U}(t)\vee
\ffilt{S_0}(t+1),
\]
from which the desired relation follows by (CI2).

Finally, to see that (GC) entails (BC), let $S=\parent{B}$ and
$A=V\without S$ for an arbitrary subset $B$ of $V$. Then the first
relation in (BC) follows directly from the global Granger-causal Markov
property. The second relation in (BC) can be derived similarly.
\end{proof}

\begin{proof}[Proof of Corollary \ref{psepgranger}]
Suppose that all paths between $A$ and $B$ are $p$-blocked given
$S$. We show that then all $B$-pointing paths between $A$ and $B$ are
$p$-blocked given $S\cup B$, which implies by the global Granger-causal
Markov property that $X_A$ is Granger-noncausal for $X_B$ with respect to
$\ffilt{A\cup B\cup S}$.

We firstly note that, in particular, every $B$-pointing path $\pi$ between
$A$ and $B$ are $p$-blocked given $S$ and, if $\pi$ does not contain
any intermediate points in $B$, also $p$-blocked given $S\cup B$.
Now, suppose that $\pi$ is a $B$-pointing path between $A$ and $B$ with
some intermediate points in $B$. Then $\pi$ can be partitioned as
$\pi=\path{\pi_1,\pi_2}$ where $\pi_1$ is a path between $A$ and some
$b\in B$ with no intermediate points in $B$. Because of the assumption,
the path $\pi_1$ is $p$-blocked given $S$ and, since it has no intermediate
points in $B$, also given $S\cup B$. It follows that all $B$-pointing
paths between $A$ and $B$ are $p$-blocked given $S\cup B$.

The other two cases $X_B\noncausal X_A\wrt{\ffilt{A\cup B\cup S}}$
and $X_A\noncorr X_B\wrt{\ffilt{A\cup B\cup S}}$ can be derived similarly.
\end{proof}

\section{p-separation in mixed graphs}
\label{psepequivalence}

The definition of $p$-separation presented in this paper
is based on paths that may be self-intersecting. This leads
to simpler conditions than in the original definition by \citet{amp2}.
The latter is formulated in terms of paths on which all intermediate
points are distinct, that is, these paths are not self-intersecting;
such paths are called trails.
According to \citet{amp2}, a trail between vertices $a$ and $b$ is
said to be $p$-active relative to $S$ if
\begin{romanlist}
\item
every $p$-collider (head-no-tail node) on $\pi$ is in $\ancestor{S}$, and
\item
every $p$-noncollider $v$ is either not in $S$ or
it has two adjacent undirected edges ($\UE v\UE$) and
$\parent{v}\without S\neq\varnothing$.
\end{romanlist}
Otherwise the trail is $p$-blocked relative to $S$.
Let $A$, $B$, and $S$ be disjoint subsets of $V$. Then $S$ $p$-separates
$A$ and $B$ if all trails between $A$ and $B$ are $p$-blocked
relative to $S$.

The following proposition shows that the two notions of
$p$-separation are equivalent.

\begin{proposition}
\label{pseparationequivalence}
Let $G=(V,E)$ be a mixed graph and $A$, $B$, $S$ disjoint subsets of $V$.
Then there exists a $p$-active trail between $A$ and $B$ relative to $S$
if and only there exists a $p$-connecting path between $A$ and $B$ given
$S$.
\end{proposition}
\begin{proof}
Suppose that $\pi$ is a trail between two vertices $a$ and $b$ that
is $p$-active relative to $S$. If all $p$-colliders on $\pi$ are
in $S$ and all $p$-noncolliders are outside $S$, then $\pi$
is also $p$-connecting given $S$. Otherwise, $\pi$ is
$p$-blocked by vertices $u_{j_1},\ldots,u_{j_r}$ on the path.
If $u_{j_i}$ is a $p$-collider then $u_{j_i}\in\ancestor{S}$
since $\pi$ is $p$-active. Hence there exists a directed
path $\tau_i=\path{u_{j_i}\DE\cdots\DE s_i}$ for some $s_i\in S$ such
that all intermediate points on $\tau_i$ are not in $S$
and we set $\sigma_i=\path{\tau_i,\rev{\tau}_i}$, where $\rev{\tau}_i$
denotes the reverse path of $\tau_i$, that is,
$\rev{\tau}_i=\path{s_i\LDE\cdots\LDE u_{j_i}}$.
On the other hand, if $u_{j_i}$ is a $p$-noncollider on $\pi$,
then the two edges adjacent to $u_{j_i}$ are undirected. Thus,
there exists $w_i\in\parent{u_{j_i}}\without S$ and we set
$\sigma_i=\path{u_{j_i}\LDE w_i\DE u_{j_i}}$. Now, let
$\pi_i$ be the subpath of $\pi$ between $u_{j_{i-1}}$ and $u_{j_i}$
with $u_{j_0}=a$ and $u_{j_{r+1}}=b$ and set
\[
\pi'=\path{\pi_0,\sigma_1,\pi_1,\sigma_2,\ldots,\pi_{r-1},\sigma_r,\pi_r}.
\]
Then all $p$-colliders on $\pi'$ are in $S$ and all $p$-noncolliders
are not in $S$, which yields that $\pi'$ is $p$-connecting given $S$.

Conversely, suppose that $\pi$ is a $p$-connecting path between
$a$ and $b$ given $S$. Let $u_{j_1}$ be the first instance of a vertex 
that occurs more than once on the path. Then $\pi$ can be partitioned as
$\pi=\path{\pi'_0,\lambda_1,\pi_1}$ such that $u_{j_1}$ is an end-point,
but not an intermediate point of $\pi'_0$ and $\pi_1$.
Noting that $\pi'_0$ is already a trail, we continue to partition
$\pi_1$ in the same way. After finitely many steps, we obtain the partition
\[
\pi=\path{\pi'_0,\lambda_1,\pi'_1,\lambda_2,\ldots,\pi'_{r-1},\lambda_r,\pi'_r}
\]
such that the subpaths $\pi'_j$ are all trails. Thus, the shortened
path $\pi'=\path{\pi'_0,\ldots,\pi'_r}$ is also a trail. We show
that $\pi'$ is a $p$-active trail relative to $S$. We firstly note
that all subtrails $\pi'_j$ are $p$-connecting and hence $p$-active.
We therefore have to show that the vertices $u_{j_i}$ satisfy the
conditions for a $p$-active trail.

Suppose that $u_{j_i}$ is
a $p$-collider that is not in $S$. Then at least one of
the edges adjacent to $u_{j_i}$ has an arrowhead at $u_{j_i}$
and we may assume that $\pi'_{i-1}$ is $u_{j_i}$-pointing
(otherwise consider the reverse path). Since $u_{j_i}\notin S$,
it must be a $p$-noncollider on $\pi$ and hence $\lambda_i$
starts with a tail at $u_{j_i}$. On the other hand, since
$u_{j_i}$ must be a $p$-noncollider on all its occurrences on $\pi$
and $\pi'_i$ does not start with a tail, the loop $\lambda_i$
cannot be a directed path (otherwise $u_{j_i}$ would not be a
$p$-collider on $\path{\lambda_i,\pi'_i}$).
Consequently there exists an intermediate point $w_i$ such that
the subpath between $u_{j_i}$ and $w_i$ is directed and $w_i$ is
a $p$-collider. It follows that $w_i\in S$ and $u_{j_i}\in\ancestor{S}$.

Next, suppose that $u_{j_i}$ is a $p$-noncollider on $\pi'$ that is in $S$.
Since $u_{j_i}$ has been a $p$-collider on $\pi$, the two edges adjacent
to $u_{j_i}$ on $\pi'$ must be undirected and $\lambda_i$ must be a
bi-pointing path. Hence $\lambda_i$ is of the from 
$\lambda_i=\path{u_{j_i}\LDE w_i,\lambda_i'}$ with $w_i\notin S$
(since $w_i$ is a $p$-noncollider and $\pi$ is $p$-connecting).
Therefore, the set $\parent{u_{j_i}}\without S$ is not empty and
$u_{j_i}$ satisfies the above condition (ii). Altogether it follows
that $\pi'$ is $p$-active relative to $S$.
\end{proof}

In a remark on our simplified version of $p$-separation,
\citet{amp2} argue that there are infinitely many possibly self-intersecting
paths in a graph as opposed to finitely many trails.
The following lemma shows that it is possible
to restrict the search for $p$-connecting paths in $G$ to a
finite number of paths, namely all paths in which no edge occurs
twice with the same orientation.

\begin{lemma}
\label{finitepaths}
Let $G=(V,E)$ be a mixed graph and suppose that $\pi$ is
a $p$-connecting path of the form $\pi=\path{\pi_1,e,\pi_2,e,\pi_3}$,
where $e$ is an oriented edge between some vertices $u$ and $v$.
Then the shortened path $\pi'=\path{\pi_1,e,\pi_3}$ is
also $p$-connecting.
\end{lemma}
\begin{proof}
Since $\pi$ is $p$-connecting, the two subpaths $\path{\pi_1,e}$ and
$\path{e,\pi_3}$ are $p$-connecting. This implies that
also $\pi'$ is $p$-connecting as every intermediate point has the
same $p$-collider/noncollider status as in the corresponding subpath.
\end{proof}

\end{appendix}
\bibliography{papers,application}
\bibliographystyle{stat}

\end{document}